\documentclass[a4paper,11pt]{article}

\usepackage{graphicx}
\usepackage{amsmath}
\usepackage{amsfonts}
\usepackage{titlesec}
\usepackage{amsthm}
\usepackage{amssymb}
\usepackage{geometry}
\usepackage{color}
\usepackage{accents}

\definecolor{dark-green}{rgb}{0.0,0.4,0.0}
\definecolor{pantone_369}{rgb}{0.4784, 0.7098, 0.0863}

\usepackage{algorithm2e}

\newcommand{\ms}{\mbox{\rm \tiny ms}}

\newcommand{\Vms}{V_{H}^{\ms}}

\newcommand{\ums}{u_H^{\ms}}

\newcommand{\dott}{(\cdot,t)}
\newcommand{\sym}{ \mbox{\tiny \rm sym} }

\newcommand{\R}{\mathbb{R}}
\newcommand{\T}{\mathcal{T}}

\newcommand{\eps}{\varepsilon}

\newcommand{\aeps}{a^\eps}
\newcommand{\ahom}{a^0}

\newcommand{\ueps}{u_\eps}
\newcommand{\hatueps}{\hat{u}_\eps}

\newcommand{\uhom}{u_0}

\newcommand{\quotes}[1]{``#1''}

\newtheorem{theorem}{Theorem}[section]

\newtheorem{lemma}[theorem]{Lemma}

\theoremstyle{definition}
\newtheorem{definition}[theorem]{Definition}
\newtheorem{remark}[theorem]{Remark}
\newcommand{\N}{\mathbb{N}}

\newcommand{\abs}[1]{{\left|{#1}\right|}}
\newcommand{\pt}{\partial_{t}}
\newcommand{\ptt}{\partial_{tt}}
\newcommand{\ps}{\partial_{s}}
\newcommand{\pss}{\partial_{ss}}
\newcommand{\pxxxx}{\partial_{xxxx}}
\newcommand{\pxx}{\partial_{xx}}

\newcommand{\pttxx}{\partial_{ttxx}}

\def\fcnL{\mathrm{L}}

\newcommand{\of}[1]{(#1)}

\def\ofOm{\of\Omega}

\newcommand{\norm}[1]{{\left\|{#1}\right\|}}

\def\sym{{\mathrm{sym}}}

\def\Ten{{\mathrm{Ten}}}
\def\Sym{{\mathrm{Sym}}}

	\newcommand{\intMeanbig}[1]{\big\langle{#1}\big\rangle}

\def\LdOm{{\Ld\ofOm}}

\newcommand{\Lp}[1]{\fcnL^{#1}}
\def\Ld{{\Lp2}}

\def\Li{{\Lp\infty}}

\def\fcnL{\mathrm{L}}

\def\Beps{\mathcal{B}^\eps}

\title{Multiscale methods for wave problems in heterogeneous media}

\begin{document}
\maketitle

\begin{center}
{\large Assyr Abdulle\footnote[1]{ANMC, Section de Math\'{e}matiques, \'{E}cole polytechnique f\'{e}d\'{e}rale de Lausanne, 1015 Lausanne, Switzerland, Assyr.Abdulle@epfl.ch} and Patrick Henning\footnote[2]{Department of Mathematics, KTH Royal Institute of Technology, SE-100 44 Stockholm, Sweden, pathe@kth.se}}\\[2em]
\end{center}

\renewcommand{\thefootnote}{\fnsymbol{footnote}}
\renewcommand{\thefootnote}{\arabic{footnote}}

\begin{abstract}
In this paper we give a survey on various multiscale methods for the numerical solution of second order hyperbolic equations in highly heterogeneous media. We concentrate on the wave equation and distinguish between two classes of applications. First we discuss numerical methods for the wave equation in heterogeneous media without scale separation. Such a setting is for instance encountered in the geosciences, where natural structures often exhibit a continuum of different scales, that all need to be resolved numerically to get meaningful approximations. Approaches tailored for these settings typically involve the construction of generalized finite element spaces, where the basis functions incorporate information about the data variations. In the second part of the paper, we discuss numerical methods for the case of structured media with scale separation. This setting is for instance encountered in engineering sciences, where materials are often artificially designed. If this is the case, the structure and the scale separation can be explicitly exploited to compute appropriate homogenized/upscaled wave models that only exhibit a single coarse scale and that can be hence solved at significantly reduced computational costs.
\end{abstract}



\section{Introduction}
\label{section:introduction}
In this article we discuss recent developments of numerical methods for the multiscale wave equation
\begin{align}
\label{eq:wave_strong}
\ptt \ueps-\nabla\cdot\left( \aeps
\nabla \ueps \right)&=F \hspace{46pt}\hbox{in}~\Omega\times]0,T[,\\
\nonumber
\ueps&=0 \hspace{49pt} \hbox{on}~\partial\Omega \times ]0,T[,\\
\nonumber
 \ueps(x,0)=g_1(x), \hspace{10pt}\pt \ueps(x,0)&=g_2(x) \hspace{30pt}\hbox{in}~\Omega,
\end{align}
where $\Omega$ is bounded domain.
We make the following minimal regularity assumptions (H0)
\begin{itemize}
\item $\Omega\subset\mathbb{R}^{d}$ is a bounded Lipschitz domain with a piecewise polygonal boundary ($d=1,2,3$);
\item the data satisfy $F \in L^2(0,T;L^2(\Omega))$, $g_1 \in H^1_0(\Omega)$ and $g_2 \in L^2(\Omega)$;
\item the matrix-valued function $\aeps$
is in $\mathcal{M}(\alpha,\beta,\Omega)$
where
\begin{eqnarray}\label{e:spectralbound}
\lefteqn{\mathcal{M}(\alpha,\beta,\Omega) :=} \\
\nonumber&\hspace{0pt}&\{ a \in [L^\infty(\Omega)]^{d\times d}_{\sym}| \hspace{5pt}
\alpha |\xi|^2 \le a(x) \xi \cdot \xi \le \beta |\xi|^2 \enspace\text{for all } \xi \in \R^d \mbox{ and almost all }x\in \Omega\}.
\end{eqnarray}
\end{itemize}
Under (H0) there exist a unique $\ueps\in C^0(0,T;H^1(\Omega))),$ 
 $\partial_{t}\ueps\in C^0(0,T;L^2(\Omega)))$ and
$\partial_{tt}\ueps \in L^2(0,T;H^{-1}(\Omega)$
such that for all $v\in H_0^1(\Omega)$  and a.e. $t > 0$ \cite[Chapter 3]{LiM72a}
\begin{eqnarray}
\label{eq:wave_weak}
\langle\partial_{tt}\ueps(\cdot,t) ,v\rangle+(\aeps \nabla \ueps(\cdot,t),\nabla v)_{L^2(\Omega)}
=(F(\cdot,t),v)_{L^2(\Omega)}.
\end{eqnarray}
The wave speed $\aeps$ is assumed to be a {\it multiscale coefficient}. By that we mean that $\aeps$ varies on a scale of order $\mathcal{O}(\eps)$, where $0<\eps \ll 1$. In general, we do not assign a particular value to $\eps$, but we only assume that it is very small parameter. However, whenever we encounter a sufficiently regular coefficient, we assume that $\| \aeps \|_{W^{1,\infty}(\Omega)}=\mathcal{O}(\eps^{-1})$ to illustrate our arguments.
For a classical numerical approximation of the problem \eqref{eq:wave_weak} we
pick a $P1$ finite element space $V_h$ with mesh size $h$ (for simplicity we consider a quasi-uniform family of triangulation of the domain $\Omega$ in simplicial elements, cf. \cite{BrS08}) and 
regard the problem: find $u_h: [0,T]\rightarrow V_h$ such that for all $v_h\in V_h$ and a.e. $t > 0$
\begin{eqnarray}
\label{eq:FEM_classical}
\langle\partial_{tt}u_h(\cdot,t),v_h\rangle+( \aeps \nabla u_h(\cdot,t),\nabla v_h)_{L^2(\Omega)}
=(F(\cdot,t),v_h)_{L^2(\Omega)}
\end{eqnarray}
and with appropriate discrete initial values. Let us from now on use the shorthand notation $L^p(H^s):=L^p(0,T;H^s(\Omega))$.
Following the arguments by Baker \cite{Bak76} leads to the standard approximation result:
\begin{eqnarray}
\label{eq:error_classical}
\|\ueps-u_h\|_{L^\infty(L^2)}\leq C (\|\ueps-\Pi_h(\ueps)\|_{L^\infty(L^2)}+\|\partial_t \ueps-\partial_t \Pi_h(\ueps)\|_{L^1(L^2)}),
\end{eqnarray}
where $\Pi_h: H_0^1(\Omega)\rightarrow V_h$ is the Ritz-projection on $V_h$, i.e., the $(a_\eps\nabla\cdot, \nabla\cdot)$-orthogonal projection.
The projection error can be further estimated by exploiting its quasi best-approximation property in $H^1(\Omega)$, an Aubin-Nitsche duality argument for the elliptic projection and the $H^1$-stability of $\Pi_h$ to obtain 
\begin{eqnarray}\label{eq:estimatefinescale}
\|\ueps-u_h\|_{L^\infty(L^2)}
&\leq & C \hspace{2pt} {h} \hspace{2pt} \|\aeps \|_{W^{1,\infty}(\Omega)} 
\left( \| \ueps \|_{L^{\infty}(H^1)} + \| \partial_t \ueps \|_{L^{1}(H^1)} \right)
\hspace{6pt} \leq \hspace{6pt} C \hspace{2pt} \frac{h}{\eps},
\end{eqnarray}
where $C=C(T)$ is independent of $\eps$. Here we assumed $\|\aeps\|_{W^{1,\infty}(\Omega)}= {\cal O}({\eps}^{-1})$ and $\| \partial_t \ueps \|_{L^{1}(H^1)}=\mathcal{O}(1)$.
At first sight, this seems not to be an optimal $L^2$ error estimate. Indeed for the Aubin-Nitsche duality argument one usually already uses an optimal first order convergence rate in $h$ (in the $H^1$ norm) for the solution of the elliptic problem. In our situation this would lead to an estimate of the type $h\|\partial_t \ueps \|_{L^1(H^2)}$.
The problem is however to bound $\|\partial_t \ueps \|_{L^1(H^2)}$ independently of $\eps$. For general initial condition classical a priori error bounds \cite{Eva10} would lead to a bound of the type $C\eps^{-2}$ \cite{AbH14c}. This is why we avoid the $H^2$ norm in the above estimate. This issue makes also the use of higher order spatial approximations useless when using classical FE or FD methods for the approximation of \eqref{eq:wave_strong}.
Notice that the standard a priori bound for $\norm{\pt\ueps}_{L^1(H^1)}$ also scales as $C/\eps$
but this term can be controlled using $G$-convergence and perturbation arguments \cite{AbH14c} 
so that the assumption $\| \partial_t \ueps \|_{L^{1}(H^1)}=\mathcal{O}(1)$ can be in fact avoided.
Hence, convergence usually requires $h<\eps$ and a high computational complexity if $\eps$ is small. In this contribution we consider two distinct situations that require different numerical strategies. For the first situation, we assume no scale separation in the highly heterogeneous tensor $\aeps$. In this case, we review numerical methods based on suitably chosen multiscale spaces. For the second situation, we assume scale separation in the tensor $\aeps$ (e.g., periodic, locally periodic, random stationary). Here we can use classical finite element spaces for the numerical approximation but the numerical methods rely on homogenized (effective) models that must be computed ``on the fly". A peculiar feature of highly oscillatory hyperbolic problems such as \eqref{eq:wave_strong} is that different homogenized models must be derived according to the time span of the desired approximation. We will consider time intervals that scale as $[0,T\eps^{-2}]$. \footnote{In the periodic homogenization setting this is the first interval of interest on which the homogenized model is not valid, e.g., for interval of the type $[0,T\eps^{-1}]$ the homogenized solution still gives adequate approximation in the $L^\infty(L^2)$ norm.} Notice that in this case the estimate \eqref{eq:estimatefinescale} for standard finite elements reads \cite{AbP15b}
\begin{eqnarray}
\label{eq:error_classical}
\|\ueps-u_h\|_{L^\infty(0,T\eps^{-2};L^2)}\leq C \frac{h}{\eps^3},
\end{eqnarray}
where again we assume a bound for $\|\partial_t \ueps \|_{L^1(H^1)}$ independent of $\eps$. Notice that for periodic problems, well-prepared initial data can be used to obtain such a bound \cite{AbH14c}.
\smallskip\\
We define the computational complexity as the size of the linear system $N$ required to be solved at each time-step $\Delta t$ of the time integrator for the wave equation. For classical FEM, in views of \eqref{eq:estimatefinescale}, we have a computational complexity of $N=\eps^{-d}$ per time step $\Delta t$
(short time interval) and a computational complexity of $N=\eps^{-3d}$ per time step $\Delta t$ in view of \eqref{eq:error_classical} (long-time interval $[0,T\eps^{-2}]$). If using an explicit method such as the popular leap-frog method for the time integration, the stability constraint $\Delta t\simeq h$  reads $\Delta t\simeq \eps$ for short-time integration and  $\Delta t\simeq \eps^3$
for long-time integration. This are indeed very severe time-step constraints due the oscillatory behavior of \eqref{eq:wave_strong}. Of course implicit methods could be used, but then the additional cost due to the linear system to solve (that again scales badly with $\eps$) constitutes a non-trivial additional cost per time-step.

We close this introduction by recalling a fundamental homogenization result (valid in both situations described above) that will be used in both classes of numerical methods described in what follows.

The basic question of classical homogenization is the following: if $(\aeps)_{\eps>0}$ represents a sequence of tensors and $(\ueps)_{\eps>0}$
the corresponding sequence of solutions to (\ref{eq:wave_weak}), does $\ueps$ converge in some sense to a limit function $\uhom$ ? Is there a (homogenized) equation for this limit function ? The hope is that due to the average limit process for $\eps \rightarrow 0$, the homogenized equation is cheap to solve. 
The abstract tool of $G$-convergence gives a general answer to that question.
\begin{definition}[$G$-convergence]\label{def-G-convergence}
A sequence $(\aeps)_{\eps >0}\subset \mathcal{M}(\alpha,\beta,\Omega)$ (i.e. with uniform spectral bounds in $\eps$) is said to be $G$-convergent to $\ahom \in \mathcal{M}(\alpha,\beta,\Omega)$ if for all $F \in H^{-1}(\Omega)$ the sequence of solutions $v^{\eps} \in H^1_0(\Omega)$ to
\begin{align*}
\int_{\Omega} \aeps \nabla v^{\eps} \cdot \nabla v = F(v) \qquad \mbox{for all } v \in H^1_0(\Omega)
\end{align*}
satisfies $v^{\eps} \rightharpoonup v^0$ weakly in $H^1_0(\Omega)$, where $v^0 \in H^1_0(\Omega)$ solves
\begin{align*}
\int_{\Omega} \ahom \nabla v^0 \cdot \nabla v = F(v) \qquad \mbox{for all } v \in H^1_0(\Omega).
\end{align*}
\end{definition}
One of the main properties of $G-$convergence is the following compactness result \cite{Spa68,Tar77}: let $(\aeps)_{\eps >0}$ be a sequence of matrices in  $\mathcal{M}(\alpha,\beta,\Omega)$, then there exists a subsequence $(a^{\eps '})_{\eps '>0}$ and a matrix $\ahom\in \mathcal{M}(\alpha,\beta,\Omega)$ such that $(a^{\eps '})_{\eps '>0}$ $G-$converges to $\ahom$. For the wave equation, we have the following result obtained in \cite[Theorem 3.2]{BFM92}:

\begin{theorem}[Homogenization of the wave equation]
\label{theorem:homogenization:wave:equation} Let assumptions (H0) be fulfilled and let the sequence of symmetric matrices $(\aeps)_{\eps >0}\subset \mathcal{M}(\alpha,\beta,\Omega)$ be $G$-convergent to some $\ahom \in \mathcal{M}(\alpha,\beta,\Omega)$. Let $\ueps \in L^{\infty}(0,T;H^1_0(\Omega))$ denote the solution to the wave equation \eqref{eq:wave_weak}. Then it holds
\begin{align*}
\ueps &\rightharpoonup \uhom \quad \mbox{weak-}\ast \mbox{ in } L^{\infty}(0,T,H^1_0(\Omega)),\\
\partial_t \ueps &\rightharpoonup \partial_t \uhom \quad \mbox{weak-}\ast \mbox{ in } L^{\infty}(0,T,L^2(\Omega))
\end{align*}
and where $\uhom \in L^2(0,T;H^1_0(\Omega))$ with $\partial_{tt}\uhom \in L^2(0,T;H^{-1}(\Omega)$ is the unique weak solution to the homogenized problem
\begin{align}
\label{homogenized-equation-weak}
\nonumber\langle \partial_{tt} \uhom \dott, v \rangle + \left( \ahom \nabla \uhom \dott, \nabla v \right)_{L^2(\Omega)} &= \left( F \dott, v \right)_{L^2(\Omega)} \qquad \mbox{for all } v \in H^1_0(\Omega) \mbox{ and } t > 0,\\
\left( \uhom( \cdot , 0 ) , v \right)_{L^2(\Omega)} &= \left( g_1, v \right)_{L^2(\Omega)} 
\qquad \hspace{21pt} \mbox{for all } v \in H^1_0(\Omega),\\
\nonumber \left( \partial_t \uhom( \cdot , 0 ) , v \right)_{L^2(\Omega)} &= \left( g_2 , v \right)_{L^2(\Omega)}
\quad \hspace{32pt} \mbox{for all } v \in H^1_0(\Omega).
\end{align}
\end{theorem}
This theorem and the compactness result stated above show that for {\it any} problem \eqref{eq:wave_weak} based on a sequence of matrices with $\aeps\in\mathcal{M}(\alpha,\beta,\Omega)$, we can extract a subsequence such that the corresponding solution of the wave problem converges to a homogenized solution. Except for special situations, e.g., locally periodic coefficients $\aeps$, i.e. tensor $\aeps(x)=a(x,\frac{x}{\eps})$ that are $\eps$-periodic on a fine scale or for random stationary tensors, it is not possible to construct $\ahom$ explicitly.

\section{Numerical methods for the wave equation in heterogeneous media without scale separation}
\label{section:nonsepration}

In this section we consider the setting that is encountered if the wave speed $\aeps$ reflects the properties of  a heterogeneous medium without scale separation. Such a setting is typically encountered in earth sciences such as geophysics or seismology. Here, the waves propagate through a medium that often lacks any kind of structure. Instead the medium consists of a variety of heterogeneously distributed materials, as for instance different rock and soil types, possibly interrupted by natural reservoirs of groundwater, petroleum or gas
(cf. the data of the Society of Petroleum Engineering, openly accessible on \texttt{http://www.spe.org/web/csp}). In such natural structures it is typically impossible to distinguish different scales of resolution and we speak about a {\it lack of scale separation}. Opposite to this, problems arising from engineering applications are often artificially designed and hence exhibit a perfect scale separation. We will discuss this case in the next section.

From now on, $\aeps$ denotes an unstructured and highly heterogeneous coefficient that lacks scale separation. In our discussion we shall also focus on the minimal possible regularity assumptions. This is an important aspect since the propagation field $\aeps$ is discontinuous in many realistic applications. Subsequently, we let $\eps$ denote a parameter that characterizes the size of the smallest length scale on which variations of $\aeps$ can be observed, e.g. we could define
$\eps := \| \aeps\|_{W^{1,\infty}(\Omega)}^{-1}$ provided that $\aeps$ is sufficiently regular. For the size of the computational domain we assume here 
that diam$(\Omega)=\mathcal{O}(1)$. This assumption will be relaxed in some situations for longtime wave propagation as described in Section \ref{sec:longtime}.
Let $X_N \subset H^1_0(\Omega)$ denote a discrete space and let $u_N: [0,T]\rightarrow X_N$ denote the corresponding semi-discrete approximation given as the solution to
\begin{eqnarray}
\label{general-discrete-problem}\langle\partial_{tt}u_N(\cdot,t),v\rangle+(\aeps \nabla u_N(\cdot,t),\nabla v)_{L^2(\Omega)}
=(F(\cdot,t),v)_{L^2(\Omega)}, \quad \mbox{for all } v \in X_N \mbox{ and a.e. } t>0,
\end{eqnarray}
and with suitably chosen discrete initial values. The question that we want to discuss is the following. Is there for any $N\in \N$ a space $X_N=X_N(\aeps)$ that only depends on $\aeps$; but not on $t$, the source term $F$ or the initial values $g_1$ and $g_2$; such that two properties are fulfilled:
\begin{itemize}
\item[{\bf(C1)}] dim$(X_N)=N$ and
\item[{\bf(C2)}] $\|\ueps-u_N\|_{L^\infty(H^1)}\leq C H$; where $H=\mathcal{O}(N^{-d})$ is a generalized mesh size and where $C=C(T)$ is independent of the variations of $\aeps$, i.e. independent of $\eps$.
\end{itemize}
The above properties would guarantee a convergence of the numerical scheme without a resolution constraint imposed by the speed of the data variations $\eps^{-1}$. Comparing this with the setting of the classical $P1$ FE space $V_h$ with corresponding solution $u_h$ given by \eqref{eq:FEM_classical} we can identify the relation dim$(V_h)=N$ with $N=\mathcal{O}(h^{-d})$, but the final error estimate 
reads $\|\ueps-u_h\|_{L^\infty(L^2)}\leq C \frac{h}{\eps}$ 
and does hence not fulfill the desired property.

Keeping these findings in mind, we next want to discuss the following problem: is there a discrete space $X_N$ that fulfills the properties {\bf(C1)} and {\bf(C2)}? In fact, the question can be answered mainly positively and we can identify four different approaches to this problem in the literature. We present them in the following in chronological order. For the sake of simplicity, the we shall also assume that the initial values are zero; i.e. $g_1=g_2=0$; and that the source term is time-independent; i.e. $F(x,t)=F(x)$. The case of general data is shortly discussed in Section \ref{section-general-initial-data}. Also note that we restrict the presentation
of the methods to the semi-discrete setting, i.e. we do not discuss a time-discretization of \eqref{general-discrete-problem} as this does typically not impose new problems. We just mention that equation \eqref{general-discrete-problem} can be discretized in time in various ways (e.g. using the framework of Newmark schemes) and that a fully discrete analysis of Approach 1 \cite{OwZ08} and Approach 4 \cite{AbH14c} can be found in the corresponding papers.

Another approach based on operator upscaling was proposed by Minkoff et al. \cite{VMK05,KoM06,VdM08}. However, we will not discuss this approach since the wave equation is considered in a different form (namely of the structure $\partial_{tt}\ueps- \aeps \triangle \ueps = F$), which does not involve the typical multiscale issues.

In the following, we only assume minimal regularity for $\aeps$, i.e. $\aeps\in [L^\infty(\Omega)]^{d\times d}_{\sym}$.

\subsection{Approach 1 - Harmonic coordinate transformations}
\label{approach-1}

Let us assume that $\Omega \subset \R^2$ is a bounded and convex domain of the class $C^2$.
The first approach devoted to the multiscale wave equation without scale separation was proposed by Owhadi and Zhang \cite{OwZ08} (as a generalization of the elliptic case considered in \cite{OwZ07}). The authors suggest to overcome the issue of the missing space regularity by applying a smoothing coordinate transformation 
$G^{\eps}=(G^{\eps}_1,G^{\eps}_2): \Omega \rightarrow \Omega,$
whose components $G^{\eps}_k$ with $k\in \{1,2\}$ are defined as the weak solutions to the elliptic boundary value problem 
\begin{align}
\label{problem-for-G-i-eps}\nonumber\nabla \cdot ( \aeps \nabla G^{\eps}_k ) &= 0 \hspace{24pt} \mbox{in } \Omega;\\
 G^{\eps}_k(x)&=x_i \hspace{20pt} \mbox{on } \partial \Omega.
\end{align} 
The transformation $G^{\eps}$ as defined above can be shown to be an automorphism over $\Omega$ (cf. \cite{AlN03}) and hence maps $\Omega$ indeed into itself.

The basic idea of the approach is to try to approximate the exact solution of \eqref{eq:wave_strong} by $v_0 \circ G^{\eps}$, where $v_0$ denotes a smooth function that is $\eps$-independent and that only exhibits slow variations. In contrast, the multiscale character and
the low regularity part of $\ueps$
are embedded in the components of the harmonic transformation $G_k^{\epsilon} \in H^1(\Omega)$. Provided that this ansatz is valid, the original multiscale problem can reinterpreted as to find the slow function $v_0(\cdot,t) \in H^1_0(\Omega) \cap H^2(\Omega)$ such that $v_0(\cdot,t) \circ G^{\eps} \in H^1_0(\Omega)$ solves \eqref{eq:wave_weak}. Since $v_0$ is a smooth single-scale function, it can be easily approximated in conventional (coarse) finite element spaces $V_H$. With that, we can define the coordinate transformed (and still low-dimensional) solution space as
\begin{align}
\label{harmonic-coordinates-space}V_H^{\eps}:=\{ v_H \circ G^{\eps}| \hspace{2pt} v_H \in V_H\}
\end{align}
and we seek $u_H^{\eps} \in (0,T; V_H^{\eps})$ with $u_H^{\eps}(\cdot,0)=\partial_t u_H^{\eps}(\cdot,0)=0$
and
\begin{eqnarray}
\label{owhadi-zhang-1}\langle\partial_{tt} u_H^{\eps}(\cdot,t) ,v\rangle+(\aeps \nabla u_H^{\eps}(\cdot,t),\nabla v)_{L^2(\Omega)}
=(F,v)_{L^2(\Omega)}, \quad \mbox{for all } v \in V_H^{\eps} \mbox{ and } t>0.
\end{eqnarray}
This approach can be rigorously justified under the following geometric assumption.
\begin{definition}[Cordes type condition]
\label{cordes-type}
Let $\sigma^{\eps}(x):= \left( \nabla G^{\eps}(x) \right)^{\hspace{-2pt}\top} \hspace{-2pt} \aeps(x) \nabla G^{\eps}(x)$,
where the columns of the matrix $\nabla G^{\eps}$ are defined by $\nabla G^{\eps}:=(\nabla G^{\eps}_1,\nabla G^{\eps}_2)$. Furthermore, the Cordes parameter is defined as
$$
\mu_{\sigma^{\eps}}:= \underset{x \in \Omega}{\mbox{ess sup}} \left( \frac{\lambda_{\max}(\sigma^{\eps}(x)) }{ \lambda_{\min}(\sigma^{\eps}(x)) } \right),
$$
where $\lambda_{\max}(\sigma^{\eps}(x))$ and $\lambda_{\min}(\sigma^{\eps}(x))$ denote the upper and lower spectral bounds of $\sigma^{\eps}(x)$.
We say that $G^{\eps}$ fulfills the Cordes type condition, if 
$\mu_{\sigma^{\eps}}<\infty$
and $\left( \mbox{Trace}[\sigma]^{-1}\right) \in L^{\infty}(\Omega)$.
\end{definition}
With this, Owhadi and Zhang \cite{OwZ08} proved the following theorem, that guarantees {\bf(C1)} and {\bf(C2)} for $V_H^{\eps}$.
\begin{theorem}\label{main-theorem-owhadi-zhang-1}
Let the Cordes type condition from Definition \ref{cordes-type} be fulfilled and assume that $g_1=g_2=\partial_t F =0$. Then we have that the coordinate-transformed solution is a regular coarse-scale function, that is $v_0:= \ueps \circ (G^{\eps})^{-1}\in L^{\infty}(0,T;H^2(\Omega))$ and
$\| v_0 \|_{L^{\infty}(H^2)} \le C \| F \|_{L^2(\Omega)}$,
i.e. the $L^{\infty}(H^2)$-norm of $v_0$ can be bounded independently of $\eps$. Furthermore, if $V_H$ denotes a conventional $P1$ finite element space as in the introduction, then the numerical approximation $u_H^{\eps}$ to \eqref{owhadi-zhang-1} (obtained in the low-dimensional multiscale space $V_H^{\eps}$) fulfills the a priori error estimate 
\begin{align*}
\| u_H^{\eps} - \ueps \|_{L^{\infty}(H^1)} + \| \partial_t u_H^{\eps} - \partial_t \ueps \|_{L^{\infty}(L^2)} \le C H \| F \|_{L^2(\Omega)},
\end{align*}
where $C$ only depends on $\Omega$ and the fulfillment of the Cordes type condition.
\end{theorem}
\begin{remark}
The error estimate in Theorem \ref{main-theorem-owhadi-zhang-1} holds with minimal regularity assumptions for $\aeps$, i.e. $\aeps\in [L^\infty(\Omega)]^{d\times d}_{\sym}$. Even though the result is only proved for $d=2$ and under the assumption that a Cordes type condition is fulfilled, both assumptions do not seem to be necessary in practice (cf. \cite{OwZ08} for a discussion and numerical experiments) and the method still performs well if the condition is not fulfilled. This is an important observation, since the validity of the Cordes-type condition can be hard to verify in practice.
\end{remark}
From the numerical perspective we can identify four steps involved in this approach.
\begin{enumerate}
\item For $k=1,\ldots,d$, solve for a numerical approximation to the components $G^{\eps}_k \in H^1(\Omega)$ of the harmonic coordinate
given by \eqref{problem-for-G-i-eps}. Since the $G^{\eps}_k$ are multiscale functions, this step involves to solve $d$ global multiscale-scale problems on a fine mesh with mesh size \quotes{$h<\eps$}.
\item When a sufficiently accurate approximation $G^{\eps,h}$ to $G^{\eps}$ is computed, we can define the basis set $\Phi_i^{\eps}:= \Phi_i \circ G^{\eps,h}$, where $\Phi_i$ denotes a nodal basis function of the coarse space $V_H$. With that it is necessary to compute the entries of the stiffness matrix $S$ and the mass matrix $M$ with $S_{ij}=(\aeps \nabla\Phi_j^{\eps}, \nabla \Phi_i^{\eps} )_{L^2(\Omega)}$ and $M_{ij}=( \Phi_j^{\eps}, \Phi_i^{\eps} )_{L^2(\Omega)}$. Note that the basis functions $\Phi_i^{\eps}$ are typically non-local, i.e. they have support in the whole domain $\Omega$.
\item By using the precomputed coarse quantities it is now possible to march in time with a favorite time-discretization, where every step only involves operations with low-dimensional matrices and vectors.
\item For new source terms $F$, the results of step 1 and 2 can be reused and it is possible to directly start with step 3.
\end{enumerate}
As we can see, Step 1 and 2 can be very costly since it involves global fine scale computations (step 1) and global fine scale quadrature rules (step 2). This can be considered as a one time overhead that pays off for a sufficiently high number of time steps or source terms. It should also be noted that the mass and stiffness matrices computed in step 2 are dense and typically not sparse. Solving a linear system that involves these matrices hence has a cubic computational complexity. This needs to be considered when deciding for a coarse mesh.
Alternativ to using conventional $P1$ finite elements for constructing the multiscale space $V_H^{\eps}$ given by \eqref{harmonic-coordinates-space}, one could also 
chose $V_H$ as the space of
weighted extended B-splines (WEB, \cite{Hoe03}). The numerical experiments in \cite{OwZ08} indicate that using B-splines can improve the performance of the method considerably. Finally, let us mention that Step 1 introduces a numerical approximation $G^{\eps,h}_{k}$ to $G^{\eps}_k$. So far, it has not yet been investigated analytically how such an approximation influences the validity of Theorem \ref{main-theorem-owhadi-zhang-1}, which assumes that $G^{\eps}_k$ is available analytically.

\subsection{Approach 2 - MsFEM using limited global information}
\label{approach-2}

A second approach that can be found in the literature is the {\it Multiscale Finite Element Method using Limited Global Information} proposed by Jiang et al. \cite{JEG10,JiE12}. This approach can be considered as a generalization of the previously discussed harmonic coordinate transformation \cite{OwZ08}, however, from the slightly shifted point of view. Jiang et al. start from the assumption that there exist $m$ known global fields $G^{\eps}_1, \dots, G^{\eps}_m$ that are available to the user. These fields might either be precomputed (e.g. coordinate transformations as in Section \ref{approach-1}) or they were inferred from measured data (for instance in the context of porous media flow). It is further assumed that there exists a smooth (unknown) function $v_0$ that is independent of the data variations, but that allows to express the exact solution $\ueps$ in terms of the global fields, i.e.
$$\ueps(x,t) \approx v_0(G^{\eps}_1(x), \dots, G^{\eps}_m(x),t).$$
Provided that this a priori knowledge is available, two approaches are proposed in \cite{JEG10} that we describe now. Let us consider a conventional coarse $P1$ finite element space denoted by $V_H$.

$\\$
{\bf Version 1.} We define the multiscale space $V_H^{\eps}$ by the functions that can be expressed as products of a global field and a finite element function, i.e. we set $G^{\eps}_0:=1$ and let
\begin{align*}
V_H^{\eps} := \{ G^{\eps}_k \hspace{2pt} v_H | \hspace{2pt} 0 \le k \le m; \hspace{4pt} v_H \in V_H \}.
\end{align*}
If $\{ \Phi_i | \hspace{2pt} 1 \le i \le n \}$ denotes the nodal basis of $V_H$, then $\{ G^{\eps}_k \hspace{2pt} \Phi_i | \hspace{2pt} 1 \le i \le n; 0 \le k \le m\hspace{4pt} \}$ denotes a basis of $V_H^{\eps}$ (with dimension $n(m+1)$). We see that the functions $G^{\eps}_k \hspace{2pt} \Phi_i$ inherit their support from $\Phi_i$, so that we obtain as set of locally supported basis functions of $V_H^{\eps}$. As before, a numerical approximation is obtained according to \eqref{general-discrete-problem}. Some analytical considerations for {\bf Version 1} are presented in \cite{JiE12}.

$\\$
{\bf Version 2.} The second version is restricted to the specific case that $d=2$ and that $m=1$, i.e. there is only one global field $G^{\eps}=G^{\eps}_1$ available and we consequently assume $\ueps(x,t) \approx v_0(G^{\eps}_1(x),t)$ for some smooth $v_0$ and for {\it continuous} $G^{\eps}_1$. In this case, Jiang et al. propose a different approach that is very similar to the classical Multiscale Finite Element Method (MsFEM, proposed by Hou and Wu \cite{HoW97}). More precisely, for every node $z_i$ of the coarse space $V_H$, we define a corresponding (multiscale) nodal basis function $\Phi_i^{\eps}$ element-wise as follows. Let $\omega_i:= \cup \{ K \in \T_H| \hspace{2pt} z_i \in K \}$ denote the nodal patch that belongs to $z_i$.
Then for every element $K \in \T_H$ of the triangulation with $K \subset \omega_i$, $\Phi_i^{\eps}$ is the solution to the elliptic problem
\begin{align*}
- \nabla \cdot \left( \aeps \nabla \Phi_i^{\eps} \right) &= 0 \hspace{23pt} \mbox{ in } K,\\
\Phi_i^{\eps} &= g^{\eps}_{i,K} \hspace{10pt} \mbox{ on } \partial K.
\end{align*}
Here, $g^{\eps}_{i,K}$ is an oscillatory boundary condition induced from the global field $G^{\eps}$, that we define in equation \eqref{msfem-glob-inf-bc} below. The classical MsFEM basis function are defined in the same way, but with the difference that the boundary condition is not oscillatory, but simply the affine condition inherited from the nodal basis, i.e. $\Phi_i^{\eps}=\Phi_i$ on $\partial K$ (cf. \cite{EfH09}). It is known that oscillatory boundary conditions typically improve the performance of the method, since the arising approximations do not suffer from so-called resonance errors (see also \cite{GrS07}). In the case of a single global field, meaningful oscillatory boundary values $g^{\eps}_{i,K}$ can be constructed in the following way. Let $K\subset \omega_i$ be a coarse element (triangle) and let $E_{K,0}$, $E_{K,1}$ and $E_{K,2}$ denote the three corresponding edges (with $\partial K = E_{K,0} \cup E_{K,1} \cup E_{K,2}$): Furthermore, we let $z_{K,0}$, $z_{K,1}$ and $z_{K,2}$ denote the three corresponding corners (nodes), where we assume (without loss of generality) that the numeration is such that $z_{K,0}=z_i$, $\partial E_{K,0} = \{ z_{K,0}, z_{K,1} \}$, $\partial E_{K,1} = \{ z_{K,1}, z_{K,2} \}$ and $\partial E_{K,2} = \{ z_{K,2}, z_{K,0} \}$. With that, we define $g^{\eps}_{i,K}$ on $\partial K$ by
\begin{align}
\label{msfem-glob-inf-bc}
g^{\eps}_{i,K}(x) :=
\begin{cases}
\frac{G^{\eps}(x) - G^{\eps}(z_{K,1})}{G^{\eps}(z_{i}) - G^{\eps}(z_{K,1})} & \mbox{if } x \in E_{K,0},\\
\hspace{30pt}0 & \mbox{if } x \in E_{K,1},\\
\frac{G^{\eps}(x) - G^{\eps}(z_{K,2})}{G^{\eps}(z_{i}) - G^{\eps}(z_{K,2})} & \mbox{if } x \in E_{K,2}.
\end{cases}
\end{align}
Assuming that the following error estimates hold
\begin{align*}
\| \ueps - v_0 \circ G^{\eps}\|_{L^{\infty}(H^1)} + \| \partial_t \left( \ueps - v_0 \circ G^{\eps} \right) \|_{L^{\infty}(L^2)} + \| \partial_{tt} \left( \ueps - v_0 \circ G^{\eps} \right) \|_{L^{2}(L^2)} \le \delta,
\end{align*}
and provided that all involved functions are sufficiently regular, it possible to derive $\delta$- and $H$-explicit and {\it$\eps$-independent} a priori error estimates for {\bf Version 2} of the method (see \cite[Theorem 3.1]{JEG10}).

$\\$
Observe that both versions do not follow the strategy suggested in \cite{OwZ08}. Comparing the {\it MsFEM using limited global information} with the approach based on a {\it harmonic coordinate transformation}, we can state two crucial differences.
\begin{enumerate}
\item Owhadi and Zhang define $V_H^{\eps}$ as the coordinate transformation of $V_H$ (hence the basis functions are concatenations $v_H \circ G^{\eps}$), whereas Jiang et al. define $V_H^{\eps}$ in a multiplicative way, which guarantees locally supported basis functions. This makes the assembling of stiffness and mass matrices significantly cheaper (in terms of quadrature costs) and leads to sparse matrix structures. Consequently, the method can be considered as cheaper, once the global fields are available. However, finding the global fields in the first place is required for both methods.
\item The harmonic coordinate transformation is supported by a rigorous numerical analysis and hence guarantees convergence of the method. On the other hand, the analysis for both versions of Approach 2 (see \cite{JiE12} for Version 1 and \cite{JEG10} for Version 2) rely on regularity and approximability assumptions that cannot be derived in general from  the transformation $G^{\eps}$ given by \eqref{problem-for-G-i-eps}.
For that reason it is not clear if there exist global fields $G^{\eps}_1$, $\cdots$, $G^{\eps}_m$ that meet the formal requirements that are necessary for the framework of Approach 2. Nevertheless, numerical experiments in \cite{JEG10,JiE12} show a good performance of Approach 2 if $G^{\eps}$ is selected as a steady state solution with $- \nabla \cdot (\aeps \nabla G^{\eps}) = F$ and a homogenous Dirichlet boundary condition. Since this choice involves the source term, the resulting space can however not be reused for different right-hand side of \eqref{eq:wave_strong}.
\end{enumerate}

\subsection{Approach 3 - Flux-transfer transformations}
\label{approach-3}

The basic drawback of Approach 1 and 2 is that a considerable one-time overhead is involved, when the harmonic coordinate transformations/the global fields are computed by a global fine scale computation. An alternative that overcomes this issue was again proposed by Owhadi and Zhang in \cite{OwZ11}, who use a localizable {\it transfer property} as an alternative to non-local harmonic coordinate transformations discussed in Section \ref{approach-1}. The transfer property says that if two fluxes $a_1 \nabla v_1$ and $a_2 \nabla v_2$ have the same divergence in $L^2$, then their potential parts (i.e. the $\nabla H^1_0(\Omega)$-parts in their Helmholtz-decomposition) can be approximated in discrete spaces with identical accuracy, if the discrete spaces $X_1$ and $X_2$ are linked through the relation $\nabla \cdot (a_1 \nabla X_1)=\nabla \cdot (a_2 \nabla X_2)$.

$\\$
In order to make this statement precise, let $\Omega\subset \R^d$ be in the following a domain with a $C^2$-boundary and let $P_{\mbox{\tiny pot}} : [L^2(\Omega)]^d \rightarrow 
\nabla H^1_0(\Omega):=
\{ \nabla v| \hspace{2pt} v \in H^1_0(\Omega) \}$ denote the $L^2$-projection onto the potential part of the Helmholtz decomposition, i.e. $P_{\mbox{\tiny pot}}(\mathbf{u})\in \nabla H^1_0(\Omega)$ fulfills
\begin{align*}
(P_{\mbox{\tiny pot}}(\mathbf{u}),\mathbf{w})_{L^2(\Omega)} = (\mathbf{u},\mathbf{w})_{L^2(\Omega)} \qquad \mbox{for all }  \mathbf{w} \in \nabla H^1_0(\Omega).
\end{align*}
With this, for any $a \in \mathcal{M}(\alpha,\beta,\Omega)$ we define the {\it $a$-flux-norm} of a function $w \in H^1_0(\Omega)$ by
\begin{align*}
\| w \|_{\mbox{\tiny-flux}} := \| P_{\mbox{\tiny pot}}(a \nabla w) \|_{L^2(\Omega)}.
\end{align*}
The flux-norm can easily seen to be equivalent to the energy-norm (with
$\alpha \| \nabla w \|_{L^2(\Omega)} \le \| w \|_{a-\mbox{\tiny flux}} \le \beta \| \nabla w \|_{L^2(\Omega)}$)
and it can be shown that it has the following remarkable transfer property (cf. \cite{OwZ11}).
\begin{lemma}\label{transfer-property-theorem}
Let $V_H\subset H^1_0(\Omega)$ be a finite dimensional subspace and let 
$$V_{H,\aeps}:= \{ \nabla \cdot (\aeps \nabla v_H)| \hspace{2pt} v_H \in V_H\}.$$
Furthermore, for $F\in L^2(\Omega)$ we let $z_{\eps}\in H^1_0(\Omega)$ and $z\in H^1_0(\Omega)$ denote the weak solutions to the following problems
\begin{align*}
\int_{\Omega} \aeps \nabla z_{\eps} \cdot \nabla v = \int_{\Omega} F \hspace{2pt} v \qquad
\mbox{and} \qquad \int_{\Omega} \nabla z \cdot \nabla v = \int_{\Omega} F \hspace{2pt} v \qquad \mbox{ for all } v \in H^1_0(\Omega).
\end{align*}
Then it holds the following transfer property in the flux norm
\begin{align}
\label{transfer-property}\inf_{v_H \in V_H} \| z - v_H \|_{\mbox{\tiny 1-flux}} = \inf_{v_{H,\aeps} \in V_{H,\aeps}} \| z_{\eps} - v_{H,\aeps} \|_{\aeps\mbox{\tiny-flux}},
\end{align}
where $\| \cdot \|_{\mbox{\tiny 1-flux}}$ denotes the flux-norm for $a=1$, i.e. $\nabla \cdot (a \nabla \cdot )= \triangle$.
\end{lemma}
To emphasize the role of the transfer property, let us denote for the rest of this subsection by $V_H$ the space of weighted extended B-splines (WEB) \cite{Hoe03} on a uniform grid with grid width $H$ (consequently the basis functions have a support on a domain with diameter $\mathcal{O}(H)$). We stress that the basis functions are smooth (hence $ \triangle \Phi_i \in L^2(\Omega)$) and locally supported and that we intrinsically assume that $V_H$ is a coarse space (i.e. the variations of $\aeps$ are not resolved). Then from the transfer property \eqref{transfer-property} and the norm equivalence we infer the $\eps$-independent estimate
\begin{align*}
\inf_{v_{H,\aeps} \in V_{H,\aeps}} \| z_{\eps} - v_{H,\aeps} \|_{H^1(\Omega)}
\le C \inf_{v_H \in V_H} \| z - v_H \|_{H^1(\Omega)} \le C H \| F \|_{L^2(\Omega)}, 
\end{align*}
where $C$ is a constant only depending on $\Omega$, $\alpha$ and $\beta$. The main advantage of Approach 3, compared to the harmonic coordinates transformation, is that the construction of $V_{H,\aeps}$ can be localized in a natural way: if $\Phi_i$ denotes a basis function of $V_H$, then the transferred basis function $\Phi_i^{\eps}$ in $V_{H,\aeps}$ is given as the solution to $\nabla \cdot (\aeps \nabla \Phi_i^{\eps}) = \triangle \Phi_i$ with $\Phi_i^{\eps}=0$ on $\partial \Omega$. Since the \quotes{source term} $\triangle\Phi_i$ is only locally supported, we expect $\Phi_i^{\eps}$ to decay to zero outside of the support of $\Phi_i^{\eps}$. This justifies that computing $\Phi_i^{\eps}$ can be localized to smaller subdomains. We shall detail this in the following.

Let $z_i \in \Omega$ denote the node associated with the basis function $\Phi_i \in V_H$ and for a sufficiently large constant $C>0$ let 
$U_i := \{ x \in \R^d| \hspace{2pt} |z_i-x| \le C H^{1/2}|\log(H)|\} $ denote an environment of $z_i$ with a diameter of order $\mathcal{O}(H^{1/2}|\log(H)|)$. Then we define the localized transferred basis function $\Psi_i^{\eps}$ as the solution to the elliptic problem
\begin{equation}\label{def-Psi-i-OZ11}
    \begin{aligned}
\frac{1}{H} \Psi_i^{\eps} - \nabla \cdot \left( \aeps \nabla \Psi_i^{\eps} \right) &= \triangle \Phi_i \hspace{30pt} \mbox{in } U_i \cap \Omega,\\
\Psi_i^{\eps} &= 0 \hspace{45pt} \mbox{on } \partial \left( U_i \cap \Omega \right).
    \end{aligned}
\end{equation}
At first glance, this problem formulation might be surprising since it involves the artificial zero-order term $\frac{1}{H} \Psi_i^{\eps}$. In fact, this term was added to speedup the decay of the corresponding Green's function. The decay of the Green's function associated with the operator 
$ \nabla \cdot \left( \aeps \nabla \cdot \right)$
is well-known to be only algebraic, whereas adding the zero-order term $\frac{1}{H} \Psi_i^{\eps}$ regularizes the operator and leads to an exponential decay of the corresponding Green's function (see also \cite{Glo11}). This exponential decay allows to restrict the computation of the transferred basis functions to computational domains $U_i$ with a size of diam$(U_i)=\mathcal{O}(H^{1/2}|\log(H)|)$. The distortion of the transfer property caused by adding the zero-order term is basically balanced with the localization error (caused by restricting the problem to $U_i$) if the constant $C$ in the definition of $U_i$ is chosen appropriately (cf. the numerical experiments in \cite{OwZ11}). Note that practically \eqref{def-Psi-i-OZ11} needs to be discretized on an additional fine mesh with mesh size \quotes{$h<\eps$}. With that, we can define the multiscale space by
\begin{align}
\label{space-owhadi-zhang-2}V_H^{\eps} := \mbox{span} \{ \Psi_i^{\eps}| \hspace{2pt} 1 \le i \le N \},
\end{align}
where $N=\mbox{dim}(V_H)=\mathcal{O}(H^{-d})$. As desired, the space is low dimensional and it has locally supported basis functions with a support of diameter $\mathcal{O}(H^{1/2}|\log(H)|)$. In particular the fact that only local problems \eqref{def-Psi-i-OZ11} have to be solved is a great advantage compared to Approach 1 and 2 which both require to first compute global fields $G^{\eps}$ by fine scale computations on the whole domain $\Omega$. The problems \eqref{def-Psi-i-OZ11} are independent and each problem individually is cheap to solve, which allows for an efficient parallel implementation of the method. Finally, the space also fulfills the properties {\bf(C1)} and {\bf(C2)} as shown by the following theorem \cite[Theorem 5.1]{OwZ11}.

\begin{theorem}\label{main-theorem-owhadi-zhang-2}
Recall that $g_1=g_2=\partial_t F =0$, the smoothness of $\Omega$ and that $V_H$ denotes the space of weighted extended B-splines on a uniform grid. If $V_H^{\eps}$ is defined according to \eqref{space-owhadi-zhang-2} and if $u_H^{\eps} \in (0,T; V_H^{\eps})$ with $u_H^{\eps}(\cdot,0)=\partial_t u_H^{\eps}(\cdot,0)=0$ solves
\begin{eqnarray*}
\langle\partial_{tt} u_H^{\eps}(\cdot,t) ,v\rangle+(\aeps \nabla u_H^{\eps}(\cdot,t),\nabla v)_{L^2(\Omega)}
=(F,v)_{L^2(\Omega)}, \quad \mbox{for all } v \in V_H^{\eps} \mbox{ and } t>0,
\end{eqnarray*}
then it holds for some $\eps$-independent constant $C>0$ that
\begin{align*}
\| u_H^{\eps} - \ueps \|_{L^{2}(H^1)} + \| \partial_t u_H^{\eps} - \partial_t \ueps \|_{L^{\infty}(L^2)} \le C H \| F \|_{L^2(\Omega)}.
\end{align*}
\end{theorem}
As a last remark we note that Theorem \ref{main-theorem-owhadi-zhang-2} does not rely on the usage of B-splines for $V_H$ and that the result is also valid for alternative choices as long as the basis function of $V_H$ are locally supported and sufficiently smooth (see \cite{OwZ11} for details).

\subsection{Approach 4 - Localized Orthogonal decomposition}
\label{approach-4-lod}

The last approach that we shall discuss was proposed in \cite{AbH14c} in the framework of the Localized Orthogonal Decomposition (LOD, cf. \cite{MaP14,HeP13}). Let $V_H$ denote a (coarse) $P1$ finite element space as in the introduction and let $P_H : H^1_0(\Omega) \rightarrow V_H$ denote the corresponding $L^2$-projection, i.e. $(P_H(v),\Phi_H)_{L^2(\Omega)}=(v,\Phi_H)_{L^2(\Omega)}$ for all $\Phi_H \in V_H$. The multiscale method in \cite{AbH14c} is derived based on the following observation (see also \cite{Pet15}):

Let $F \in L^2(\Omega)$ and let $z_{\eps} \in H^1_0(\Omega)$ denote the solution to the elliptic multiscale problem
\begin{align*}
\int_{\Omega} \aeps \nabla z_{\eps} \cdot \nabla v = \int_{\Omega} F v \qquad \mbox{for all } v \in H^1_0(\Omega).
\end{align*}
Then the $L^2$-projection of $z_{\eps}$ (i.e. the $L^2$-best approximation of $z_{\eps}$ in $V_H$) can be characterized as the unique solution to the Petrov-Galerkin problem
\begin{align}
\label{LOD-motivation}\int_{\Omega} \aeps \nabla P_H(z_{\eps}) \cdot \nabla (\mbox{Id}+Q)(v_H) = \int_{\Omega} F \hspace{2pt} (\mbox{Id}+Q)(v_H) \qquad \mbox{for all } v_H \in V_H
\end{align}
and where $-Q: V_H \rightarrow \mbox{ker} P_{H}:=W$ denotes the $(\aeps \nabla \cdot, \nabla \cdot )_{L^2(\Omega)}$-orthogonal projection into the kernel of the $L^2$-projection. 
Indeed, defining $$V_H^{\eps}:=\{ v_H +Q(v_H)| \hspace{2pt} v_H \in V_H\},$$ we observe that $H^1_0(\Omega)=V_H\oplus W$  with $~V_H  \perp W$ with respect to the $(\cdot , \cdot)_{L^2(\Omega)}$ scalar product while $H^1_0(\Omega)=V_H^{\eps}\oplus W$ with $~V_H  \perp W$ with respect to the $(\aeps \nabla \cdot, \nabla \cdot )_{L^2(\Omega)}$ scalar product.

 As problem \eqref{LOD-motivation} is a finite dimensional problem (involving only the degrees of freedom from the coarse space $V_H$), it can be solved cheaply, provided that an approximation to the operator $Q$ is available. In particular, we have the regularity- and $\eps$-independent error estimate
$$
\|Êz_{\eps} - P_H(z_{\eps}) \|_{L^2(\Omega)} \le C H \| z_{\eps} \|_{H^1(\Omega)} \le C H \| F \|_{L^2(\Omega)}.
$$

Motivated by the above considerations, a low-dimensional discrete multiscale space is constructed as follows. Let $V_h$ denote a conventional $P1$ finite element space with fine mesh size \quotes{$h<\eps$} (so that all variations of $\aeps$ are resolved) and such that $V_H \subset V_h$. We denote the kernel of $L^2$-projection $P_H$ by $W_h:=\{ w_h \in V_h| \hspace{2pt} P_H(w_h) =0\}$.
\begin{enumerate}
\item In order to approximate the Ritz-projection $Q_h : V_h \rightarrow W_h$ given by 
$$(\aeps \nabla Q_h(v), \nabla w_h )_{L^2(\Omega)}= -(\aeps \nabla v, \nabla w_h )_{L^2(\Omega)}$$
for all $w_h \in W_h$ in an efficient way, we make an affine decomposition of the right hand side into $(\aeps \nabla v, \nabla w_h )_{L^2(\Omega)}=\sum_{K\in \T_H}(\aeps \nabla v, \nabla w_h )_{L^2(K)}$ and solve for each of the localized sources $(\aeps \nabla v, \nabla w_h )_{L^2(K)}$ individually. Since the corresponding solutions can be shown to exhibit an exponential decay outside of $K$, we can replace the computational domain $\Omega$ by a small environment $U_k(K)$ of $K$. Here, $U_k(K)$ is defined iteratively by
\begin{equation}\label{def-patch-U-k}
    \begin{aligned}
      U_0(K) & := K, \\
      U_k(K) & := \cup\{T\in \T_H\;\vert\; T\cap U_{k-1}(K)\neq\emptyset\}\quad k=1,2,\ldots,
    \end{aligned}
\end{equation}
i.e. the patch $U_k(K)$ consists of the coarse element $K$ and $k$-layers of coarse elements around it. With that, we define 
$W_h(U_k(K)):=\{ v_h \in W_h| \hspace{2pt}
v_h=0 \enspace \mbox{in } \Omega \setminus
U_k(K) \}$
and solve for $Q_h^K(v_H) \in W_h(U_k(K))$ with
\begin{align}
\label{corrector-problems}\int_{U_k(K)} \aeps \nabla Q_h^K(v_H) \cdot \nabla w_h =- \int_{K}  \aeps \nabla v_H \cdot \nabla w_h \qquad \mbox{for all } w \in W_h(U_k(K))
\end{align}
and we set the global approximation $Q_{h,k}$ of $Q_h$ to
$$
Q_{h,k}(v_H):= \sum_{K\in \T_H} Q_h^K(v_H).
$$
The stationary local problems \eqref{corrector-problems} are fully independent from each other and can be hence solved in parallel. Furthermore, for small $k$, they are of small size and hence only involve a low computational complexity. Note that problem \eqref{corrector-problems} only has to be solved for coarse basis functions that have a support in $K$. Details on how they can be solved practically are provided in \cite{AbH14c} and \cite{EHM16}.
\item 
With the precomputed operator $Q_{h,k}$ we define the multiscale space by
\begin{align}
\label{LOD-space}\Vms := \{ v_H +Q_{h,k}(v_H)| \hspace{4pt} v_H \in V_H\}.
\end{align}
The space is low-dimensional and if $k$ is large enough so that $Q_{h,k}=Q_{h}$, we obtain the $(\aeps \nabla \cdot, \nabla \cdot )_{L^2(\Omega)}$-orthogonal decomposition $V_h = \Vms \oplus W_h$. As $k$ is typically only a small number, we speak about a {\it Localized Orthogonal Decomposition} (LOD).
\item March in time with a favorite time-discretization, where every step only involves operations in $\Vms$ with dim$\Vms$=dim$V_H$.
\item For new source terms $F$, the results can be reused.
\end{enumerate} 
The Localized Orthogonal Decomposition approach can be shown to fulfill the following error estimates \cite{AbH14c}.
\begin{theorem}\label{a-priori-error-estimates-LOD}
Let the localization parameter $k\in \N$ be chosen such that $k\simeq|\log(H)|$ and
let $\Vms$ denote the corresponding LOD multiscale space. Then, if $\ums \in (0,T; \Vms)$ with $\ums(\cdot,0)=\partial_t \ums(\cdot,0)=0$ solves
\begin{eqnarray*}
\langle\partial_{tt} \ums(\cdot,t) ,v\rangle+(\aeps \nabla \ums(\cdot,t),\nabla v)_{L^2(\Omega)}
=(F,v)_{L^2(\Omega)}, \quad \mbox{for all } v \in \Vms \mbox{ and } t>0,
\end{eqnarray*}
there exists a constant $C$ that only depends on the time $T$, on $\Omega$ and on the mesh regularity, such that the following error estimates hold true
\begin{align*}
\|  \ueps - \ums \|_{L^{\infty}(L^2)} \le C H^2 \| F \|_{L^2(\Omega)} + 
e_h^{\mbox{\tiny$(1)$}}
\end{align*}
and
\begin{align*}
\| \partial_t\ueps - \partial_t \ums \|_{L^{\infty}(L^2)} 
+ \| \ueps -  \ums \|_{L^{\infty}(H^1)} \le C H + e_h^{\mbox{\tiny$(2)$}},
\end{align*}
where $e_h^{\mbox{\tiny$(1)$}} := \| \ueps - \Pi_h(\ueps) \|_{L^{\infty}(L^2)} + \| \partial_t \ueps - \Pi_h(\partial_t \ueps) \|_{L^1(L^2)}$ and $
e_{\mbox{\tiny \rm disc}}^{\mbox{\tiny$(2)$}}(h) := \| \partial_t \ueps - \Pi_h(\partial_t \ueps) \|_{L^{\infty}(L^2)} 
+ \| \partial_{tt} \ueps - \Pi_h(\partial_{tt} \ueps) \|_{L^1(L^2)} +  \|\ueps - \Pi_h(\ueps) \|_{L^{\infty}(H^1)}
$ denote the fine scale discretization errors, with $\Pi_h$ being the Ritz-projection on $V_h$.
\end{theorem}
From Theorem \ref{a-priori-error-estimates-LOD} we see that the LOD allows for optimal error estimates in $L^{\infty}(L^2)$, $L^1(L^2)$ and $L^{\infty}(H^1)$ without additional regularity assumptions. Furthermore, we see that the precomputations necessary to construct $\Vms$ involves small elliptic problems in computational domains of size $\mathcal{O}(H |\log(H)|)$. Since the problems can be solved in parallel, the method can be implemented efficiently. As for Approach 3, the fact that the computations for Approach 4 can be localized is a considerable advantage compared to Approaches 1 and 2. In terms of the size of the local domains, we see that Approach 4 involves subdomains that are by the factor $\sqrt{H}$ smaller than the ones necessary for Approach 3. On the down side, the local problems \eqref{corrector-problems} in Approach 4 are saddle point problems which involves the computation and inversion of the Schur complement matrix associated with the constraint \quotes{$P_H(w)=0$} (cf. \cite{EHM16} for details). Approach 3 on the contrary only involves unconstrained local problems. As a final difference, the discrete spaces required for Approach 3 need to be of higher order and smooth (e.g. weighted extended B-splines), whereas Approach 4 can be implemented using conventional $P1$ Lagrange finite element spaces.

$\\$
Recently it has been shown that Approach 4 also intrinsically relaxes the CFL condition on adaptive meshes \cite{PeS16}, which is very significant for corresponding time-discretizations. Generalizations of the approach to the Helmholtz equation in the context of high frequency wave propagation are given in \cite{BGP15,GaP15,Pet16}.

\subsection{The case of general initial values: G-convergence and perturbation arguments}
\label{section-general-initial-data}

In this section we shall discuss the case of general initial values, i.e. $g_1$, $g_2$, $\partial_t F\neq 0$. First, we easily observe that all the Approaches 1-4 can be straightforwardly modified to fit this case. Typically the discrete initial values are picked as the $L^2$- or Ritz-projections onto the multiscale space. Why is it therefore necessary to discuss this case independently? To understand the issue, note that if $g_1=g_2=\partial_t F= 0$, then all higher order time derivatives of $\ueps$ in $t=0$ will vanish as well, i.e. we have (as far as it exists) $\partial_t^{j} \ueps(\cdot,0)=0$ for all $j\ge 0$. However, if we consider the general case, we figure out that for $j\ge 2$ the time-derivatives in $t=0$ are linked through 
$
\partial_t^j \ueps(\cdot,0) = \partial_t^{j-2} F(\cdot , 0) + \nabla \cdot ( \aeps \nabla \partial_t^{j-2} \ueps(\cdot,0)).
$
From that expression we can see that the \quotes{higher order initial values} can be rapidly oscillating with a high amplitude of order $\mathcal{O}(\eps^{1-j})$ for $j\ge 2$. Unfortunately, these terms will just pop up on the right hand side of the a priori error estimates for Approach 1-4, i.e. in the Theorems \ref{main-theorem-owhadi-zhang-1},  \ref{main-theorem-owhadi-zhang-2} and \ref{a-priori-error-estimates-LOD} (see also \cite{AbH14c,JiE12,JEG10,OwZ08,OwZ11}) and the convergences rates can again depend on $H/\eps$, what we just wanted to avoid. 

In fact, there is currently no positive result on the question if {\bf(C1)} and {\bf(C2)} can be still fulfilled simultaneously for general initial values. 
However, a positive answer has been recently given when we replace the $L^{\infty}(H^1)$ in {\bf(C2)} by the 
$L^{\infty}(L^2)$ norm.

As for elliptic multiscale problems, it is easy to verify that the $L^p$-norm of any space derivative of $\ueps$ that is higher than $1$, will explode with decreasing $\eps$, i.e. we can say that for any $s>1$ we have in general that $\| \ueps \|_{L^1(H^s)} \overset{\eps \rightarrow 0}{\longrightarrow} \infty$.
It is tempting to assume that the $L^1(H^s)$-norms (for $s>1$) are the only norms that should be avoided. However, unfortunately this is not the case and the statement can be \quotes{often} generalized to $\| \ueps \|_{W^{m,2}(H^s)} \rightarrow \infty$ for $\eps \rightarrow 0$, whenever $m+s>1$ (cf. \cite{CiD99}). The vague quantification \quotes{often} is detailed by the following theorem (cf. \cite{AbH14c}), which directly links the problem to the choice of initial values.

\begin{theorem}[Time-regularity and regularity estimates]\label{proposition-time-regularity}
Let $\ueps$ denote the solution to the wave equation \eqref{eq:wave_weak} and assume that $F \in W^{m,2}(0,T;L^2(\Omega))$ for some $m\in \mathbb{N}$. Recursively we define the generalized initial values $w_j^{\eps}\simeq \partial_t^j \ueps(\cdot,0)$ by
\begin{align}
\label{higher-order-initial-values}w_0^{\eps}:=g_1, \qquad w_1^{\eps}:=g_2, \qquad w_j^{\eps} := \partial_t^{j-2} F(\cdot , 0) + \nabla \cdot ( \aeps \nabla w_{j-2}^{\eps}) \quad \mbox{for } j=2,3,\cdots,m+1. 
\end{align}
If $w_j^{\eps} \in H^1_0(\Omega)$ for $0\le j \le m$ and $w_{m+1}^{\eps}\in L^2(\Omega)$, we have
$$
\partial_t^{m} \ueps \in L^{\infty}(0,T;H^1_0(\Omega)); \enspace  \partial_t^{m+1} \ueps  \in L^{\infty}(0,T;L^2(\Omega))
\enspace \mbox{and} \enspace  \partial_t^{m+2} \ueps \in L^{2}(0,T;H^{-1}(\Omega))
$$
and it holds the (optimal) regularity estimate
\begin{eqnarray}
\label{reg-estimates}
\nonumber\lefteqn{\| \partial_t^{m} \ueps \|_{L^{\infty}(0,T;H^1(\Omega))} + \| \partial_t^{m+1} \ueps \|_{L^{\infty}(0,T;L^2(\Omega))}}\\
&\le&C \hspace{2pt} \left( \| F \|_{W^{m,2}(0,T;L^2(\Omega))} + \| w_m^{\eps} \|_{H^1(\Omega)} +  \| w_{m+1}^{\eps} \|_{L^2(\Omega)} \right),
\end{eqnarray}
where $C>0$ is a constant that only depends on $T$ and $\Omega$.
\end{theorem}
From the theorem we see that we can only hope for $\eps$-independent bounds for e.g. $\| \partial_t \ueps \|_{L^{\infty}(H^1)}$ and $\| \partial_{tt} \ueps \|_{L^{\infty}(L^2)}$, if the initial values $g_1$ and $g_2$ are picked such that they cancel out the variations of $\aeps$. As considered in the previous sections, a trivial case is the case with $g_1=g_2=0$ and $F \in L^2(\Omega)$ being constant in time, for which we conclude $\partial_t^{m} \ueps \in L^{\infty}(0,T;H^1_0(\Omega))$ for all $m\in \mathbb{N}$ and $\| \partial_t^{m} \ueps \|_{L^{\infty}(0,T;H^1(\Omega))} \le C \hspace{2pt} \| F \|_{L^2(\Omega)}$. This guarantees that $\eps$-oscillations cannot enter through time derivatives. 

In general, we observe that the right hand side of the regularity estimate \eqref{reg-estimates} can be only bounded independently of $\eps$ (i.e. independent of the speed of the variations of $\aeps$) if $g_1$, $g_2$ and $\partial_t F$ are trivial or if $g_1=g_1^\eps$, $g_2=g_2^\eps$ and $F(\cdot,0)=F^{\eps}(\cdot,0)$ are {\it $\eps$-dependent} and such that they interact adequately with the variations of $\aeps$. This assumption can be hard to verify in practice. Also observe that {\it if not} $g_1=g_2=\partial_t F(\cdot,0)= 0$, the assumptions of Theorem \ref{proposition-time-regularity} for $m\ge 1$ can typically only be fulfilled if $\aeps$ admits a higher order regularity such as $\aeps \in W^{1,\infty}(\Omega)$. This can be problematic in realistic applications, where the propagation field $\aeps$ is often discontinuous. The missing smoothness of $\aeps$ can hence be a further issue in addition to the multiscale character.

An argument to overcome the problems that arise form general initial values is presented in \cite{AbH14c}. Here it is proposed to slightly perturb the initial value $g_1$ in such a way, that the effect is almost \quotes{invisible} in $L^2$; i.e. $\| g_1 - g_1^{\eps} \|_{L^2}\le \delta \ll 1$; but such that the variations of the perturbed initial value $g_1^{\eps}$ interact adequately with the variations of $\aeps$; i.e. $\| \nabla \cdot ( \aeps \nabla g_1^{\eps})\|_{L^2(\Omega)} = \mathcal{O}(1)$. By stability arguments, the solution $\hatueps$  to the perturbed initial value
is close to the original solution $\ueps$, in the sense that $\| \hatueps - \ueps \|_{L^{\infty}(L^2)} \le \| g_1 - g_1^{\eps} \|_{L^2}\le \delta$. Since the second order initial value $\partial_{tt} \ueps(\cdot,0) = F(\cdot , 0) + \nabla \cdot ( \aeps \nabla g_1^{\eps})$ is harmless (it does no longer blow up with decreasing $\eps$), the arguments for the trivial case apply again.

The construction of $g_1^{\eps}$ is obtained via $G$-convergence (see Definition \ref{def-G-convergence}). For instance, for the LOD-approach (Section \ref{approach-4-lod}) the following result has been proved in \cite{AbH14c}.
\begin{theorem}
With the assumptions and the notations of Theorem \ref{a-priori-error-estimates-LOD}, we assume that $\aeps$ is $G$-convergent with $G$-limit $\ahom$ 
and we let $g_1^{\eps} \in H^1_0(\Omega)$ denote the solution to
\begin{align}\label{def-f-eps}
\int_{\Omega} \aeps \nabla g_1^{\eps} \cdot \nabla v = \int_{\Omega} \ahom \nabla g_1 \cdot \nabla v \qquad \mbox{for all } v \in H^1_0(\Omega),
\end{align}
i.e. $g_1^{\eps}$ is constructed such that its homogenized limit coincides with $g_1$. If $g_2 \in H^1_0(\Omega)$, $\partial_t F \in L^2(0,T,L^2(\Omega))$ and $\nabla \cdot (\ahom \nabla g_1) + F(\cdot,0)\in L^2(\Omega)$ then it holds
\begin{align*}
\lim_{h\rightarrow 0} \|  \ueps - \ums \|_{L^{\infty}(L^2)} \le C H \left( \| F \|_{W^{1,2}(L^2)} + \| g_1 \|_{H^1} + \| g_2 \|_{H^1} 
+ \| \nabla \cdot (\ahom \nabla g_1) 
\|_{L^2} \right)
+ \| g_1^\eps - g_1 \|_{L^2},
\end{align*}
where $C=C(T)$.\end{theorem}
Since $\| g_1^\eps - g_1 \|_{L^2}\rightarrow 0$ for $\eps \rightarrow 0$ by the definition of $G$-convergence, we can assume that  $\| g_1^\eps - g_1 \|_{L^2}\le C H$ and hence the total observable convergence rate will by linear in $H$, i.e. $\lim_{h\rightarrow 0} \|  \ueps - \ums \|_{L^{\infty}(L^2)} \lesssim H$.
Note that $g_1^{\eps}$ in \eqref{def-f-eps} is never computed, but is just a tool useful for the analysis.
We close this section by mentioning that the above arguments do no longer apply for $L^{\infty}(H^1)$-errors.

\section{Numerical methods for the wave equation in heterogeneous media with scale separation}
\label{section:separation}
When the oscillatory tensor $\aeps$ in the wave equation \eqref{eq:wave_strong} exhibit scale separation, numerical homogenization methods with much lower computational cost than the methods described in Section \ref{section:nonsepration} can be constructed. A typical situation of scale separation for the oscillatory coefficient $\aeps$ is for example a locally periodic structure, i.e., when
$\aeps$ is of the form $\aeps(x)=a(x,x/\eps)=a(x,y)$, where $a(x,y)$ is $Y$-periodic in $y$. Here $Y$ is a  unit cell, e.g., $Y=(0,1)^d$. Another situation is the case of a tensor $\aeps$ modeled by 
a random field $\aeps=a(x/\eps,\omega)$, where the analog of the periodicity in this case is the stationarity of the field. In this review we will focus on locally periodic tensor and refer to \cite{FGP07} for a recent account on the theory and numerics for wave in random media.

The theoretical framework for the development of numerical methods in locally periodic media is that of 
$G$-convergence. Indeed in locally periodic media, it is known that for each $x\in\Omega$ 
the whole sequence of symmetric oscillatory tensors $(\aeps)_{\eps >0}$ G-converges to a unique  effective tensor $\ahom(x)$. Hence the whole family of solutions $\ueps$ of the wave equation \eqref{eq:wave_strong}
converges to the solution of an effective wave equation \eqref{homogenized-equation-weak}. It can be shown that
$\ahom$ is again in $\mathcal{M}(\alpha,\beta,\Omega)$ and
\begin{equation}
\label{equ:a_0}
\ahom_{ij}(x) = \frac{1}{|Y|}\int_Y e_i^T \hspace{2pt} a(x,y)(\nabla\chi_j(x,y) + e_j)dy,
\end{equation}
where $\chi^j(x,y),~j=1,\ldots d$  are the solutions in  $H_{per}^1(Y)$ of  
so-called cell problems
\begin{equation}
\label{equ:cell}
\int_Y a(x,y)(\nabla\chi_j(x,y)+e_j)\cdot\nabla w(y) dy=0~~\hspace{30pt} \mbox{for all } w \in H_{per}^1(Y),
\end{equation}
where $e_j$ are the vectors of the canonical basis of $\mathbb{R}^d$. Observe that the solution of problem
\eqref{equ:cell} is unique up to a constant that needs to be fixed.
Except for the case when $\aeps$ is periodic, the full computation of the map $x\rightarrow \ahom(x)$ is not possible as it relies on infinitely many PDE solutions. Hence a numerical homogenization scheme must rely on a finite set of well chosen homogenized tensors $\ahom(x_j)_{j=1}^M$. 

\subsection{Effective model and numerical homogenization method for short-time wave propagation}

A number of numerical methods based on homogenization theory and the asymptotic expansions \cite{BLP78} have been recently proposed. We mention the method given in \cite{DoC09} that  applies to problems with uniformly periodic oscillatory tensors that are assumed to be symmetric with respect to the center of the periodic cell. Furthermore the macroscopic computational domain is assumed to be the union of an integer number of periodic cells. Numerical homogenization methods for the wave equation have also appeared in the geoscience community.
In \cite{GCM10} a numerical homogenization method based on asymptotic expansion \cite{BLP78} is derived for elastic waves. Unlike the FE-HMM algorithm presented below that couples macro and micro scales in a global numerical scheme, the procedure in \cite{GCM10} is sequential and consists in first finding an effective tensor (analytical expression are used for 1d periodic problems, or spectral element method are used for 2d problems) and then solving the effective wave equation. It should be noted that a filtering technique to compute effective parameters for highly oscillatory non-periodic wave equations is discussed. While no rigorous theoretical foundation of such procedure seems to be available, this technique is shown to give good results for some problems. Such filtering techniques could also be used in a pre-processing step for the numerical homogenization algorithm described below. Recently, a finite difference method (FDM) \cite{EHR11}
and a finite element method \cite{AbG11} for wave problems in highly oscillatory media have been proposed in the framework of heterogeneous multiscale methods (HMM). These HMM methods are general algorithms to approximate numerically the homogenized solution in case of scale separation in the oscillatory tensor $\aeps$ that we will describe in more detail.
We discuss here finite element or finite difference algorithms but note that HMM methods have also be coupled  with spectral element methods \cite{AbE07}. Such a coupling could readily be implemented for the methods described below.

\subsubsection{Approach  1 - Finite-element numerical homogenization method}
\label{sec:FE_HMM}
In this section we describe the finite element heterogeneous multiscale method (FE-HMM) proposed in \cite{AbG11} that can be seen as a general numerical homogenization method.
We pick a standard macroscopic finite element space $V_H$ made of piecewise polynomial on each macro element $K$
of a partition of the computational domain $\Omega=\cup_{K\in{\cal T}_H}K$. We next define a sampling domain $K_{\delta}$ (of size $\delta$ comparable to $\eps$) within each macro element $K$. A micro finite element space with a triangulation that resolves the fine scale $\eps$ is defined in each sampling domain. We consider then the following problem: 
find $u_H: [0,T]\rightarrow V_H$ such that
\begin{eqnarray}
\label{equ:wave_HMM}
& &(\partial_{tt} u_H(\cdot,t),v_H)
+B_H(u_H(\cdot,t),v_H)=(F(\cdot,t),v_H) \hspace{30pt} \mbox{for all }  v_H\in V_H,
\end{eqnarray}
with appropriate projection of the true initial conditions, where
\begin{eqnarray}
\label{equ:B_H}
B_H(u_H,v_H):=  \sum_{K \in{\cal T}_H} \frac{\abs{K}}{\abs{K_\delta}}\int_{K_\delta} 
a^\varepsilon(x)\nabla u_{K}^h \cdot\nabla v_{K}^h dx,
\end{eqnarray}
and $u_{K}^h$ (respectively $v_{K}^h$) are solutions of the following micro problems: for $K \in {\cal T}_H$ find $(u_{K}^h-u_H)\in V_h(K_{\delta})$ such that
\begin{equation}
\label{equ:cell-discrete}
\int_{K_{\delta}}\
~a^\varepsilon(x)\nabla u_{K}^h \cdot \nabla z^h dx=0
\hspace{30pt} \mbox{for all } z^h\in V_h(K_{\delta}).
\end{equation}

\begin{remark}
The above method is a particular case of a more general algorithm proposed in \cite{AbG11}. Indeed, $V_H$ can be replace by $V_H^\ell$ in which $v_H |_K$ is a either a polynomial of  total degree $\ell$ if $K$ is a simplex or  a polynomial
of degree at most $\ell$ in each variable if $K$ is a parallelogram. We then need $j=1,\ldots,J$ integration points 
and sampling domains $K_{\delta_j}=x_{K_j}+\delta I,$ where $I=(-1/2,1/2)^d$
and $\delta_j$ the size of the sampling domain is such that $\eps \le \delta_j \ll H$.
The bilinear form \eqref{equ:B_H} becomes
\begin{equation}
\label{equ:wave_HMM_b}
B_H(u_H,v_H)=  \sum_{K\in{\cal T}_H} \sum_{j=1}^J
\frac{\omega_{K_j}}{|K_{\delta_j}|} \int_{K_{\delta_j} }
a^\varepsilon(x)\nabla u_{K_j}^h \cdot\nabla v_{K_j}^h dx,
\end{equation}
where $u_{K_j}^h$ are solutions of \eqref{equ:cell-discrete} in ${K_{\delta_j} }$ such that $(u_{K_j}^h-u_{H,\hbox{lin}})\in V_h^q(K_{\delta_j})$, where $u_{H,\hbox{lin}}$ is a piecewise linear approximation of $u_H$ in $K$ around $x_{K_j}$
and $V_h^q(K_{\delta_j})$ is a $q$-th order finite element space.
\end{remark}
A fully discrete analysis of the method \eqref{equ:wave_HMM} has been obtained in \cite{AbG11}. Let the assumptions (H0) be fulfilled and consider a sequence $(\aeps)_{\eps >0}\subset \mathcal{M}(\alpha,\beta,\Omega)$ that $G$-converges to $\ahom \in \mathcal{M}(\alpha,\beta,\Omega)$. 
Then for the approximation of the solution $\uhom$ of the corresponding wave equation  \eqref{homogenized-equation-weak} by the FE-HMM solution
\eqref{equ:wave_HMM} with the general modified bilinear form \eqref{equ:wave_HMM_b} it holds
\begin{theorem}
Under suitable regularity of the solutions $u_{K_j}^h$ to the micro-cell problems and of the homogenized tensor $\ahom$ the error
$e_H=\uhom-u_H$ satisfies
\begin{eqnarray}
\label{equ:427}
 \|\pt e_H\|_{L^{\infty}(0,T;L^2(\Omega))}+\|e_H\|_{L^{\infty}(0,T;H^1(\Omega))}
&\leq&
C\left(H^l + \Big( \frac h\varepsilon \Big)^{2q}+err_{mod}\right),\\
\label{equ:428}
\|e_H\|_{L^{\infty}(0,T;L^2(\Omega))} &\leq&
C\left( H^{l+1} +\Big( \frac h\varepsilon \Big)^{2q} +err_{mod}\right),
\end{eqnarray}
where $C$ is independent of $H,h,\eps$ but depends on $T$.
\end{theorem}
The term $err_{mod}$ stands for the  modeling error and quantifies how well the micro averaging procedure captures the homogenized tensor at a given quadrature point $x_{K_j}$. This error depends on the size of the sampling domain $K_{\delta_j}$, the boundary condition used to solve \eqref{equ:cell-discrete}, the nature of the scale separation in $\aeps$ (e.g., local periodicity, random stationary, etc.). For example, for locally periodic coefficients, if the measure of $K_{\delta_j}$ is an integer number of cubes of size $\eps$ and $V_h^q(K_{\delta_j})$ is a subspace of a Sobolev space of periodic functions, we have $err_{mod}=0.$ If $V_h^q(K_{\delta_j})$ is a subspace of $H_0^1(K_{\delta_j})$ and $\delta>\eps$ arbitrary, then $err_{mod}=C(\delta+\eps/\delta)$ \cite{AbG11}. Further results and recent improvement of this bounds also for random stationary problems are reported in \cite{GlH14}.
\smallskip\\
We conclude this subsection by discussing the computational complexity per time-step $\Delta t$ for the FE-HMM. For simplicity of the discussion, let $\ell=q=1$.
Assuming that $N_{macro}=1/H$, $H$ is the diameter of the macroscopic FE, the size of the system of ODEs that need to be solved at each time step $\Delta t$ is $N=H^{-d}$ for piecewise linear macro FEs. We denote this computational complexity per time step as $\mathrm{cost}(\Delta t,H^{-d})$. Notice that the stability constraint if using an explicit time-integrator reads
$\Delta t\simeq H \gg\eps$. Hence the macroscopic time-step can be chosen independently of $\eps$ in sharp contrast with the  resolved FEM for the original highly oscillatory wave equation \eqref{eq:wave_strong}.
Next we discuss the offline cost for the FE-HMM: this is the cost involved in solving the micro problems in each sampling domain $K_{\delta_j}$. We note that this computation has only to be performed once at the beginning of the macroscopic time integration, unless the highly oscillatory tensor $\aeps$ is time dependent. 
Assume that we use $N_{micro}$ 
elements in each space dimension for the discretization of the sampling domains $K_{\delta_j}$, then 
$h={\delta}/{N_{micro}}$ and $h/\eps=\delta/(\eps N_{micro})={\cal O} (1/{N_{micro}})$ as $\delta\simeq \eps$.
Hence the microscopic degrees of freedom read $M_{micro}=(N_{micro})^d$ and the total offline cost is given by the solution of ${\cal O}(H^{-d})$ linear systems of size $(N_{micro})^d$. As micro and macro error must be balanced, in view of the estimates \eqref{equ:427}, \eqref{equ:428} we can choose $N_{micro}=\sqrt{N_{macro}}$
to guarantee a macroscopic linear convergence rate in the energy norm and  $N_{micro}=N_{macro}$ for a macroscopic quadratic convergence rate in the $L^2(L^\infty)$  norm.

\subsubsection{Approach  2 - Finite-difference numerical homogenization method}
\label{sec:FD_HMM}
In this section we describe a finite difference heterogeneous multiscale method (FD-HMM) proposed in \cite{EHR11}  that can be seen as a general finite difference numerical homogenization method.
The FD-HMM follows the methodology first given in \cite{AbE03} and is also related to the HMM for first oder hyperbolic problems given in \cite{CES05}.
Consider the homogenized equation \eqref{homogenized-equation-weak} in strong form
(where we set $F\equiv 0$ for simplicity)
\begin{align}
\label{eq:wave_strong_hom}
\ptt \uhom &=\nabla\cdot 
\left( \ahom \nabla \uhom \right)
\hspace{94pt}\hbox{in}~\Omega\times]0,T[,\\
\nonumber
\uhom &=0 \hspace{148pt} \hbox{on}~\partial\Omega \times ]0,T[,\\
\nonumber
 \uhom(x,0)&=g_1(x), \hspace{10pt}\pt \uhom(x,0)=g_2(x) \hspace{30pt}\hbox{in}~\Omega.
\end{align}
Consider a spatial grid $x_{i_1,\ldots,i_d}~,i_1\ldots,i_d=1,\ldots,N_{macro}$ with mesh size $H$ for $\Omega$. The first step is the macroscopic discretization of the flux formulation for the effective equation. For simplicity of exposition we set $d=2$
and $(i_1,i_2)=(i,j)$.
We seek the evolution of a function 
$U(t)=\{U_{i,j}(t)\}$ such that
\begin{eqnarray}
\label{equ:wave_HMM_fd_mac}
& & \frac{d^2}{dt^2}U_{i,j}(t)=
\frac{1}{H}(J_{{i+\frac{1}{2},j}}^H(t)-J_{{i-\frac{1}{2},j}}^H(t))+\frac{1}{H}(J_{{i,j+\frac{1}{2}}}^H(t)-J_{{i,j-\frac{1}{2}}}^H(t)),
\end{eqnarray}
where
 $J_{{i\pm\frac{1}{2},j\pm \frac{1}{2}}}^H(t)\simeq 
 \ahom(x_{i\pm\frac{1}{2},j\pm\frac{1}{2}}) \nabla u^0(x_{i\pm\frac{1}{2},j\pm\frac{1}{2}},t)
$.
The approximation of $\nabla \uhom(\cdot,t)$ at the points $x_{i\pm\frac{1}{2},j\pm\frac{1}{2}}$ is denoted by 
$P_{i\pm\frac{1}{2},j\pm\frac{1}{2}}^H(t)$. For example 
\begin{eqnarray}
\label{equ:wave_HMM_flux_mac}
P_{i,j+\frac{1}{2}}^H=
\begin{pmatrix}
\frac{(U_{i+1,j+1}-U_{i+1,j})-(U_{i-1,j+1}-U_{i-1,j})}{4H}  \\
\frac{U_{i,j+1}-U_{i,j}}{H} \\
 \end{pmatrix},
\end{eqnarray}
and similarly for the other gradients $P_{i\pm\frac{1}{2},j\pm\frac{1}{2}}^H$.

The evaluation of $J_{{i\pm\frac{1}{2},j\pm \frac{1}{2}}}^H(t)$ is obtained by solving a micro problem that relies, as for the FE-HMM, on the original multiscale data. For each domain $K_{mic}=K_{\delta}(x_{i\pm\frac{1}{2},j\pm\frac{1}{2}}):=x_{i\pm\frac{1}{2},j\pm\frac{1}{2}} + \delta I$ 
 we solve for the oscillatory wave equation \eqref{eq:wave_strong} 

\begin{eqnarray}
\label{eq:wave_micro}
& &\pss u_{mic}(x,s)=\nabla\cdot\left( \aeps(x)
\nabla u_{mic}(x,s) \right)~~\hbox{in}~K_{\delta}(x_{i\pm\frac{1}{2},j\pm\frac{1}{2}})\times[0,\tau],\\
\nonumber
& & u_{mic}(x,0)=P_{i\pm\frac{1}{2},j\pm\frac{1}{2}}^H(t)\cdot x,~~ \ps u_{mic}(x,0)=0,~~\hbox{on}~~K_{\delta}(x_{i\pm\frac{1}{2},j\pm\frac{1}{2}})\times\{s=0\},\\
\nonumber
& & u_{mic}(x,s)-u_{mic}(x,0)~~\hbox{is}~K_{\delta}(x_{i\pm\frac{1}{2},j\pm\frac{1}{2}})~\hbox{-periodic}.
\nonumber
\end{eqnarray}
The computation of the macroscopic flux $J_{_{i\pm\frac{1}{2},j\pm \frac{1}{2}}}^H$ is then obtained by computing the space time average of $\aeps\nabla u_{mic}$ over the domains $K_{\delta}(x_{i\pm\frac{1}{2},j\pm\frac{1}{2}})\times[0,\tau]$. The computation of this macroscopic flux can be further improved by computing suitable weighted averages of the microscopic flux. This can reduce the error coming from the artificial boundary conditions (the so called resonance or modeling error) and the error arising from the time averaging of the oscillatory functions
(we refer to \cite{EHR11} for details). 
\begin{remark}
Of course in practice \eqref{eq:wave_micro} needs to be discretized numerically. In \cite{EHR11}, the same numerical method as for the macro problem \eqref{equ:wave_HMM_fd_mac} is used. Precisely, we consider
a spatial grid $x_{k,l}~,k,l=1,\ldots,N_{micro}$ with meshsize $h$ for each domain $K_{\delta}$ and
seek the evolution of a micro function 
$u=\{u_{k,l}\}$ such that
\begin{eqnarray}
\label{equ:wave_HMM_fd_mic}
& &\frac{d^2}{ds^2}u_{k,l}(s)=
\frac{1}{h}(J_{{k+\frac{1}{2},l}}^h(s)-J_{{k-\frac{1}{2},l}}^h(s))+\frac{1}{h}(J_{{k,l+\frac{1}{2}}}^h(s)-J_{{k,l-\frac{1}{2}}}^h(s)),
\end{eqnarray}
together with the initial and boundary condition of \eqref{eq:wave_micro}.
Here  $J_{{k\pm\frac{1}{2},l\pm \frac{1}{2}}}^h=\aeps(x_{{k\pm\frac{1}{2},l\pm \frac{1}{2}}})p_{{k\pm\frac{1}{2},l\pm \frac{1}{2}}}$, where $p_{{k\pm\frac{1}{2},l\pm \frac{1}{2}}}$ is the micro gradient that can be computed with a similar formula as in \eqref{equ:wave_HMM_flux_mac} with $u_{k,l}$ instead of $U_{i,j}$ and $h$ instead of $H$. The weighted average of the discrete micro gradients is then used to compute the macro fluxes $J_{{i\pm\frac{1}{2},j\pm \frac{1}{2}}}^H$.
We note that  $h<\eps$ is required to compute these micro problems.
\end{remark}

In \cite{EHR11}, the method presented above has been analyzed for periodic problems in one dimension assuming $H/\eps\in\mathbb{N}$
and the following error estimate has been derived 
\begin{eqnarray}
\label{err_FD-HMM_l2}
\sup_{0\leq t_n\leq T}|U_i^n-u^0(x_i,t_n)|\leq C(H^2+err_{mic}+err_{mod})\quad\hbox{for all $x_i\in\Omega$},
\end{eqnarray}
where $C$ depends on $T$. The analysis presented in \cite{EHR11} assumes an exact micro flux $\aeps\nabla u_{mic}$. A generalized analysis was recently provided in \cite{ArR16}.
Taking into account the discretization error introduced by \eqref{equ:wave_HMM_fd_mic} would allow further to quantify $err_{mic}.$ Notice that the modeling error can be quantified for periodic problems and depends on the choice of the weighted average used to compute the micro flux averages.
\smallskip\\
Similarly as for the FE-HMM, assuming a macroscopic mesh size $H$ the
size of the system of ODEs that need to be solved at each time step $\Delta t$ is $N=H^{-d}$. However at first sight, the evaluation of the fluxes need to be performed at each macro time step, representing an additional cost of
$\tau/(\delta t)$ steps of a micro time integrator solving an ODE of size $M_{micro}=(N_{micro})^d$ corresponding to \eqref{eq:wave_micro} ($N_{micro}$ is independent of $\eps$ and corresponds to the number of grid points in each spatial dimension of the mesh used for the sampling domain $K_{\delta}$). But as the computation of the macroscopic flux $J^H$ depends linearly on the macroscopic gradient $P^H$ through \eqref{eq:wave_micro}, we observe that $J^H$ can be expanded in a linear combination of precomputed macroscopic fluxes $J^H$ solution of \eqref{eq:wave_micro} with initial conditions $u_{mic}(x,0)=e_i\cdot x,$ where $e_i,~i=1,\ldots,d$ is the canonical basi of $\mathbb{R}^d$. This is of course only valid if the oscillatory tensor $\aeps$ is time-indepent, in which case 
the solutions  of the micro problems in the FD-HMM can be considered as an offline cost.

\subsection{Effective model and numerical homogenization method for  long-time wave propagation}
\label{sec:longtime}
It has been observed for some time \cite{SaS91,FCN02} that while classical homogenization describes well wave propagation in heterogeneous medium for short time, i.e.,
\begin{equation}
\label{hom_short}
\norm{\ueps-\uhom}_{\Li(0,T;\LdOm)} \leq C\eps, 
\end{equation}
dispersive effects accumulate for longer time in the oscillatory wave $\ueps$ that are not captured by the homogenized wave $\uhom$. Precisely, the estimate \eqref{hom_short} is no longer valid for time interval
$[0,T\eps^{-2}].$ For periodic problems, using formal Bloch waves techniques, a higher order effective model has been derived in \cite{SaS91} capable of capturing the dispersive effects of the true wave over longer time intervals. Based on this effective model, a generalization of the FD-HMM method \eqref{sec:FD_HMM} has been derived in \cite{EHR12}. 
However, as the effective model derived in \cite{SaS91,FCN02} is ill-posed, the related numerical method needs to be regularization. The whole numerical scheme is not trivial to analyze, and while in \cite{ArR14} the flux error of the generalized FD-HMM scheme has been analyzed for one-dimensional periodic problems, a full analysis of the numerical scheme is yet to be done.
Recently another dispersive limit has been derived in \cite{Lam11b}, again for periodic problems, in the form of a Boussinesq type equation and the first rigorous error analysis for time intervals $[0,T\eps^{-2}]$ for the error
between $\ueps$ and the solution of the Boussinesq equation has been derived.
Multi-dimensional problems have been analyzed in
  \cite{DLS15} using Bloch waves techniques. The approach in \cite{Lam11b} has  been generalized in \cite{AbP15b}, where a whole family of effective equations has been derived and rigorous error estimates of the corresponding solutions towards the highly oscillatory wave $u^\eps$ have been established. Generalization for
multi-dimensional problems have appeared in  \cite{AbP16}. A generalization of the FE-HMM, called FE-HMM-L valid over long-time has been proposed in \cite{AGS13,AGS14}. In \cite{AbP15b} a rigorous analysis of the FE-HMM-L has been given. In the following we briefly describe various effective models valid over long-time intervals and discuss corresponding numerical schemes.

\subsubsection{A family of effective equations for the wave equation over long-time}
Consider \eqref{eq:wave_strong} with  $\aeps=a(x/\eps)=a(y)$ is $Y$-periodic ($Y$ is a unit cube, e.g., $(0,1)^d$).
In this situation the homogenized model is given according to  \eqref{homogenized-equation-weak} by
\begin{align}
\label{homogenized-equation-weak_l}
 \partial_{tt} \uhom  + \ahom_{ij}\partial^2_{ij}\uhom = F,
\end{align}
that we write here in a slightly different form  using Einstein summation rule. 
The constant  homogenized tensor $\ahom$ is given by \eqref{equ:a_0} and the initial conditions
for \eqref{homogenized-equation-weak_l} are those given in \eqref{homogenized-equation-weak}. For time intervals of length $[0,T\eps^{-1}]$
it is still possible to show that \cite{AbP15b}
\begin{equation}
\label{hom_teps}
\norm{\ueps-\uhom}_{\Li(0,T\eps^{-1};\LdOm)} \leq C\eps, 
\end{equation}
where $C$ depends on $T$ but not on $\eps$. Hence the numerical methods described in Sections
\ref{sec:FE_HMM},\ref{sec:FD_HMM} are still a good approximation of $u^\eps$ in the $L^\infty(L^2)$. For example for the FE-HMM we obtain
 \begin{equation}
\label{FE-HMM_teps}
\norm{\ueps-u_H}_{\Li(0,T\eps^{-1};\LdOm)} \leq C\left( H^{l+1} +\Big( \frac h\varepsilon \Big)^{2q}\right),
\end{equation}
where for periodic problem, the modeling error error in \eqref{equ:428} vanishes
if one uses periodic FE spaces for the micro solver with an integer number of oscillatory periods in each direction of the sampling domains.
The model \eqref{homogenized-equation-weak_l} however does not longer capture effective behavior of the true
oscillatory wave $\ueps$ for time intervals $[0,T\eps^{-2}]$ \cite{SaS91,FCN02}. In turn, an estimate such as 
\eqref{hom_teps} is no longer valid.
Building on  \cite{Lam11b,DLS15} and following \cite{AbP16} we consider for a domain $\Omega$ that is a union of cells of volume $\eps^d|Y|$ the following effective equation
\begin{equation}	\label{eq:effectiveEquation}
\begin{array}{ll}
\ptt\tilde u - \ahom_{ij}\partial^2_{ij}\tilde u 	+ \eps^2\big(a^2_{ijkl}\partial^4_{ijkl}\tilde u - b^0_{ij}\partial^2_{ij}\ptt\tilde u\big) = F,
\end{array}
\end{equation}
such that $x\mapsto\tilde u(t,x)$ is $\Omega\text{-periodic}$
with initial conditions given in \eqref{homogenized-equation-weak}. Notice that this is an effective model for the equation  \eqref{eq:wave_strong} with periodic boundary conditions on $\Omega$. We also note that most of the results below are valid on arbitrary large domain $\Omega$ with appropriate changes \cite{AbP15b,DLS15}. Indeed two related model for wave propagation over long-time can be considered: in the first model the wave propagates on long-time interval through a large spatial domain without hitting its boundary, in the second model, the spatial domain is ${\cal O}(1)$ and the wave enters and leaves many times the domain thanks to the periodic boundary conditions. One crucial difference in the analysis is that in the former case, one needs to control the Poincar\'e inequality arising from the large diameter of $\Omega$. This can be done by choosing an appropriate weak norm in space. We refer to \cite{AbP15b,DLS15} for details.
In the above approximation, the coefficients $\ahom_{ij}$ are the homogenized coefficients that appear in \eqref{homogenized-equation-weak_l}. 
The coefficients
$b^0_{ij}$ are the entries of a matrix $b^0\in \hbox{Ten}^2(\mathbb{R}^d),$ while the coefficients $a^2_{ijkl}$ are the entries of a tensor of order four $a^2\in \hbox{Ten}^4(\mathbb{R}^d)$, where $\hbox{Ten}^n(\mathbb{R}^d)$ is the space of tensors of order $n$ in $\mathbb{R}^d$. We also denote by $\Sym^n(\R^d)$ the subspace of
$\Ten^n(\R^d)$ of symmetric tensors.
If we assume that 
 \begin{equation}
\label{tensor_ass}
b^0\in\Sym^2(\R^d),~b^0 \eta \cdot \eta \geq 0~~\forall\eta\in\R^d;~~ a^2_{ijkl}=a^2_{lkji},~a^2(\eta\eta^T) : (\eta\eta^T) \geq0~~\forall\eta\in\R^d,
\end{equation}
 then (recall that $\ahom$ is elliptic and bounded) the problem \eqref{eq:effectiveEquation} is well posed
(see \cite{AbP16}). The approximation property of a solution to \eqref{eq:effectiveEquation} is summarized in the following theorem proved in \cite{AbP16}. 
Recall that $\chi_i$ is a solution to Problem \eqref{equ:cell}.
To analyze  the solution to \eqref{eq:effectiveEquation} over long-time intervals, further families of cell-problems need to be considered. First define a so-called adaptation operator of the form
\begin{equation}	\label{eq:ansatzAsymptoticExpansionWave}
\Beps\tilde u(t,x) = \tilde u\big(t,x\big) + \eps \chi_i(y)\partial_{i}\tilde u(t,x) +\eps^2 \theta_{ij}(y)\partial^2_{{ij}}\tilde u(t,x) +\ldots,
\end{equation}
where we note that an expansion of up to order four in $\eps$ is needed, that we skip here for simplicity.
Plugging this approximation into \eqref{eq:effectiveEquation} and 
separating the different
powers of $\eps$ 
allows to define appropriate equations for $\chi_i,\theta_{ij}$. First we see
that 
$\chi_i$, for $i=1,\ldots,d$, is a solution to Problem \eqref{equ:cell} (order $\eps^{-1}$), second at order $\eps^0$ we obtain that
$\theta_{ij}$, for $i,j=1,\ldots d$, are solutions of 
\begin{equation}
\label{equ:cell_2}
\int_Y a(x,y)\nabla\theta_j(y)\cdot\nabla w(y) dy=G_{ij}(w)~~\forall w \in H_{per}^1(Y),
\end{equation}
where $G_{ij}(w)=
S^2_{ij}\big\{- ( a(x,\cdot) e_i \chi_j , \nabla w)_{L^2(Y)} + ( a(x,\cdot) (\nabla\chi_j + e_j) - \ahom(x) e_j,  e_i w )_{L^2(Y)}  \big\}
$ and  $S^n_{ij}(\cdot)$ denotes the symmetrization  operator of a given tensor. Continuing this procedure up to order $\eps^{2}$ leads to the following conditions on the tensors
$b^0$ and $a^2$ that guarantee an accurate approximation of the solution of \eqref{eq:wave_strong} by the solution of \eqref{eq:effectiveEquation} over long-time. We summarize results obtained in
\cite{AbP15b,AbP16}	.				
\begin{theorem}
Under appropriate regularity assumptions on the data $\aeps,g_1,g_2$ and
if $b^0$ and $a^2$ satisfy the relation
\begin{equation}	\label{eq:constraintforeffectiveequation}
S^4_{ijkl}\big\{ a^2_{ijkl}-\ahom_{ij}b^0_{kl}\big\}
= 
S^4_{ijkl}\big\{
\intMeanbig{  a_{jk}\chi_l\chi_i}_Y
	-\intMeanbig{a\nabla\theta_{ji}\cdot \nabla\theta_{kl}}_Y
		-\ahom_{jk}\intMeanbig{\chi_l\chi_i}_Y
	\big\},
\end{equation}
where $\chi_l$ is any solution of \eqref{equ:cell} and $\theta_{kl}$ is any corresponding solution of \eqref{equ:cell_2}.
Then, the following error estimate holds
\begin{equation}	
\label{eq:thmeffectivenessTildeUb_estimate}
\norm{\ueps-\tilde u}_{\Li(0,T\eps^{-2};\LdOm)}
	 \leq C \eps 	
	 \end{equation}
where $C=C(T)$ is independent of $\eps$.
\end{theorem}

For one-dimensional problems, among the class of effective equations, there is a very simple representative that reads \cite{Lam11b,AbP15b}
\begin{equation}	
\label{eff:simple}
\ptt \tilde u - \ahom \pxx \tilde u - \eps^2 b^0 \pttxx \tilde u = F,
\end{equation}
where $b^0=\frac{1}{|Y|}\int_Y \chi^2 dy$ and $\chi$ is the solution of the cell problem \eqref{equ:cell}.

\subsubsection{Numerical homogenization methods for the wave equations over long-time}
We start by describing the generalization of the FE-HMM introduced in \cite{AGS13,AGS14}.
Recall that the FE-HMM is given by \eqref{equ:wave_HMM}. Next we replace the  $L^2$ scalar product
$(\partial_{tt} u^H,v^H)$ with $(\partial_{tt} u^H,v^H)_Q$ where 
$$
(u_H,v_H)_Q:=(u_H,w_H)+\sum_{K\in{\cal T}_H}\frac{|K|}{|K_{\delta}|}
\int_{K_{\delta}} (u_{K}^h-v_{H})
(v_{K}^h-v_{H}) dx,
$$
where $u_{K}^h$ (respectively $v_{K}^h$)  are the micro functions already used in \eqref{equ:B_H}. The FE-HMM-L method is then defined by: find $u_H: [0,T\eps^{-2}]\rightarrow V_H$ such that
\begin{eqnarray}
\label{equ:wave_HMM-L}
& &(\partial_{tt} u_H,v_H)_Q
+B_H(u_H,v_H)=(F,v_H)~~~\forall v_H\in V_H.
\end{eqnarray}
We observe that the cost of this method is similar to the original FE-HMM as the modification of the $L^2$ scalar product relies on micro functions that need anyway to be computed to assmble $B_H(\cdot,\cdot)$.
As for the FE-HMM we can replace $B_H(\cdot,\cdot)$ by the more general bilinear form
\eqref{equ:wave_HMM_b} and likewise $(u_H,v_H)_Q$ can be replaced by
$$
(u_H,v_H)_Q=(u_H,w_H)+\sum_{K\in{\cal T}_H} \sum_{j=1}^J
\frac{\omega_{K_j}}{|K_{\delta_j}|} \int_{K_{\delta_j} }
( u_{K_j}^h-u_{H,\hbox{lin}}))(v_{K_j}^h -v_{H,\hbox{lin}}))dx.
$$
For one-dimensional problems, an error analysis over long-time of the FE-HMM-L has been given in  \cite{AbP15b}.
Under suitable regularity assumptions  it has been shown that
\begin{eqnarray}
\label{equ:FE-HMM-L_err}
\|u^\eps-u_H\|_{L^\infty(0,\eps^{-2}T;L^2(\Omega))}\leq C\Big( \epsilon + \big({h}/{\eps^2}\big)^2  + H^{\ell+1}/\eps^2+H^\ell/\eps \Big)
\end{eqnarray}
where $C$ is independent of $\eps$. In fact it can be seen that for one-dimensional problems, the FE-HMM-L approximates the effective equation \eqref{eff:simple}. We emphasize that in sharp contrast with the full fine-scale approximation with error estimate \eqref{eq:error_classical}, we can take higher order FEs for the FE-HMM-L and do not require that $H<\eps$, while higher order FEs in \eqref{eq:error_classical} would result in additional negative power of $\eps$ as discussed in the introduction. Furthermore even in the regime $H<\eps$ the computational cost of the FE-HMM-L is much smaller than the cost of applying a classical FEM applied to \eqref{eq:wave_strong} as described below. For multi-dimensional problems, however, the fourth order term in the effective equation \eqref{eq:effectiveEquation} does not vanish and the FE-HMM-L cannot approximate the full effective equation. In \cite{AbP16}, an efficient numerical algorithm 
to approximate the effective coefficients in \eqref{eq:effectiveEquation} based on the solution of $d + \binom{d+1}{2}$ cell problems is given. The resolution of the effective equation relies on a fast Fourier transform algorithm. 

A generalization of the FD-HMM proposed in \cite{EHR11} has also been derived in order to capture long-time dispersive effects \cite{EHR12,ArR14}. This method has only been investigated for one-dimensional problems and we briefly describe the method for this case. Similarly to \eqref{equ:wave_HMM_fd_mac} we assume $F=0$ and consider a FD macroscopic flux formulation for the evaluation of $U=\{U_{i}\}$ that reads
\begin{eqnarray}
\label{equ:wave_HMM_fd_mac_long}
& & \frac{d^2}{dt^2}U_{i}=
\frac{1}{24H}(J_{i-3/2}^H-27J_{i-1/2}^H+27J_{i+1/2}^H-J_{i+3/2}^H).
\end{eqnarray}
For a given smooth function, the finite difference approximation based on the above scheme ensure the fourth order approximation of its derivative. The evaluation of $J_{i\pm 3/2}^H, J_{i\pm 1/2}^H$ relies on micro problems similar to \eqref{eq:wave_micro} but with higher order initial value, i.e., $u_{mic}(x,0)=q(x)$, where $q(x)$ is a cubic polynomial obtained by interpolating the current macroscopic solution at four points around the micro sampling domains. Some care is required to define the actual $q(x)$ used for computation as it is based on a subtle post-processing of the actual cubic interpolation polynomial  \cite{EHR12,ArR14}. Finally, the actual computation of $J_{i\pm 3/2}^H, J_{i\pm 1/2}^H$ is again based on a weighted average of $\aeps\nabla u_{mic}$ over the domains $K_{\delta}(x_{i\pm\frac{1}{2}})\times[0,\tau]$. A more accurate procedure as for short-time approximation is also required for this averaging procedure.
At a given point $x$, the computed flux  $J_{x}^H$ flux is shown to approximate the flux 
$$
\hat F=\ahom \partial_x \hat u+\eps^2 b^0\partial_{xxx}\hat u
$$
that is the flux of the ill-posed effective equation derived in \cite{SaS91} that reads for $F=0$ (compare with \eqref{eff:simple})
\begin{equation}	
\label{eff:simple_ill}
\ptt \hat u - \ahom \pxx \hat u - \eps^2 b^0 \pxxxx \hat u =0.
\end{equation}
\smallskip\\
Finally, we discuss
the computational complexity of the FE-HMM-L for one-dimensional problems, as this is the only method for which an a priori error analysis is available that in turn gives an estimate on the size of the spatial mesh size 
needed for long-time integration. In view of the estimate \eqref{FE-HMM_teps} for $\ell=1$ the size of the linear system to be solved per time-step $\Delta t$ is
$N=\eps^{-1}$. The CFL constraint reads here $\Delta t\simeq \eps$. This is significantly less expensive than the cost per time step over long-time intervals by
a classical FE solver, for which $N=\eps^{-3}$ and the CFL constrain reads $\Delta t\simeq \eps^3$.
Furthermore, for the FE-HMM-L higher order macro solvers can be used and hence $H$ can be chosen larger than $\eps$ (see the error estimate \eqref{equ:FE-HMM-L_err}). Here again, as for the FE-HMM, the solution of the micro problems is a one shot offline cost that is negligible in view of the macroscopic cost for each time step $\Delta t$ over a long-time time interval $[0,T\eps^{-2}]$.
In contrast, higher order FEM for a direct fine scale solver applied to \eqref{eq:wave_strong} does not improve the time cost as already discussed in Section \ref{section:introduction}.

$\\$
{\bf Acknowledgements.}
This work is partially supported by the Swiss National Foundation.

\def\cprime{$'$}

\end{document}